
\documentstyle{amsppt}
\magnification = \magstep 1
\pageheight{7 in}
\refstyle{A}
\input xypic

\topmatter
\title
Del Pezzo Surfaces over Dedekind schemes
\endtitle
\author
Alessio Corti
\endauthor
\address Department of Mathematics, University of Chicago, 5734 S.
University Avenue, Chicago IL 60637 \endaddress
\email corti\@math.uchicago.edu \endemail
\toc
\specialhead 1. Introduction \endspecialhead
\specialhead 2. Integral models of cubic surfaces \endspecialhead
\specialhead 3. Geometric integral models of Del Pezzo surfaces of degree
$d\geq 3$ \endspecialhead
\specialhead 4. Del Pezzo surfaces of degree $2$ \endspecialhead
\specialhead 5. Examples of local rigidity \endspecialhead
\specialhead References \endspecialhead
\endtoc
\endtopmatter

\newpage

\document

\head 1. Introduction \endhead

\item{a.} Discussion of the problem

\item{b.} Statement and discussion of the results

\item{c.} Open questions

\item{d.} Structure of the paper

\item{e.} Acknowledgments
\medskip

{\bf a. Discussion of the problem.}

I will begin recalling the definition of Del Pezzo surfaces, then
state the main question studied in this paper. The results are theorems
1.10, 1.15 and 1.18 below. The main ideas and techniques of proof are
introduced in \S 2.

\definition{1.1 Definition} Let $K$ be a field. A Del Pezzo surface is a
Gorenstein surface $X_K$ over $K$, with ample anticanonical bundle.
\enddefinition

The most important discrete invariant of a Del Pezzo surface is its degree,
the selfintersection number of the canonical class. It is well known that
$1 \leq \operatorname{deg} \leq 9$. Smooth Del Pezzo surfaces have been
known classically, for the singular ones see \cite{De} \cite{HW} \cite{R6}.
I will now introduce some of the notation and conventions in use throughout
this paper.

\definition{1.2 Basic setup}
We fix a Dedekind scheme $S$, with
fraction field $K$.
We denote $\eta=\operatorname{Spec} K$ the generic point of $S$.
In what follows all
varieties, schemes, morphisms, etc. are tacitly assumed to be defined over
$S$, unless otherwise indicated. For a scheme $Z$, $Z_\eta$ is the generic
fiber. A {\it birational map} is always assumed
to be {\it biregular} when restricted to the generic fibers.
\enddefinition

We wish to study the following problem:

\definition{1.3 Problem} Given a smooth Del Pezzo surface $X_K$ over $K$,
find a ``nice'' integral model of $X_K$, i.e., a nice scheme $X$ over
$S$ such that $X_\eta = X_K$.
\enddefinition

It is important to realize that we do not allow base
change by a finite field extension $L\supset K$, since most questions
about
Del Pezzo surfaces, both in
arithmetic and in geometry, become uninteresting after such a base change
(for a discussion of moduli of Del Pezzo surfaces, see \cite{I}).

We obtain definite results only when $S$ is a normal complex curve
(the ``geometric case''), but, at least for cubic Del Pezzo surfaces,
we have conjectural results for (nearly) arbitrary $S$ (2.11).

As already hinted, there are two sources of motivations to study 1.3,
geometry and arithmetic. In the arithmetic context, $K$ is a number field
and 1.3 was raised first, to my knowledge, in \cite{MT}. My own motivation to
consider 1.3 comes from the birational classification theory
of 3-folds, namely the study of Mori fiber spaces of dimension 3.
These are fibrations whose general fibers have ample anticanonical class,
the higher dimensional analogue of (old) minimal rational and ruled
surfaces. Recall that a variety is ${\Bbb Q}$-factorial if every Weil
divisor is ${\Bbb Q}$-Cartier. Terminal singularities are defined in 2.5.

\definition{1.4 Definition}
A Mori fiber space is a projective morphism
$f: X \to T$, where $X$ is a projective variety with
${\Bbb Q}$-factorial terminal singularities, and $f: X \to T$ is an
extremal contraction of fibering type. This means that:

1.4.1 $f$ has connected fibers, $T$ is normal and $\operatorname{dim} (T)
< \operatorname{dim} (X)$ (in general, $f$ is not equidimensional).

1.4.2 $-K_X$ is $f$-ample and the rank of the relative Picard group is
$1$: $\rho (X/T) = \rho(X) - \rho (T) =1$.
\enddefinition

The requirements that $X$ be ${\Bbb Q}$-factorial and $\rho (X/T)=1$ are
crucial (see e.g. the proof of the N\"other-Fano inequalities in \cite{Co}).

When $X$ is 3-dimensional, Mori fiber spaces come in three kinds,
according to the dimension of the base space $T$:

{\bf 1.4.2.1} $\operatorname{dim}(T)=0$. In this case $X$ is a {\it ${\Bbb
Q}$-Fano 3-fold} with Picard number $\rho (X)=1$.

{\bf 1.4.2.2} $\operatorname{dim}(T)=1$. $X \to T$ is a (flat) {\it Del Pezzo
fibration}, in the sense that the generic fiber is a smooth Del Pezzo
surface. This is the case of most interest to us.

{\bf 1.4.2.3} $\operatorname{dim}(T)=2$. $X \to T$ is a {\it conic bundle},
in the sense that the generic fiber is a smooth conic. $f$ is not
necessarily flat, but all fibers are curves.

The very notion of Mori fiber spaces, and the special place they occupy in
classification theory, arise from, and are best illustrated by, the
statement of the minimal model theorem for 3-folds \cite{Mo3},
which also suggests a first (quite crude) answer to 1.3:

\proclaim{1.5 Minimal Model Theorem for 3-folds}
Let $X$ be a smooth projective
3-fold over the field of complex numbers.
Then a sequence of elementary birational operations $X \dasharrow X^1
\dasharrow \cdots \dasharrow X^m = X^\prime$, called
divisorial contractions and flips,
creates a projective 3-fold $X^\prime$, with
${\Bbb Q}$-factorial terminal singularities, such that
precisely one of the following alternatives occurs:

1.5.1 $X^\prime$ is a minimal model, that is, $K_{X^\prime}$
is nef.

1.5.2 $X^\prime$ is a Mori fiber space, that
is, there exists a contraction of an extremal ray $f : X^\prime
\to T$, to a lower dimensional normal variety $T$, satisfying the
conditions of 1.4

The birational equivalence class of $X$ determines which of the
alternatives occurs.
\endproclaim

The realization that, in a Mori fiber space, it is the
extremal contraction $f: X \to T$ that should be considered the relevant
structure, rather than the variety $X$ itself, is the fruit of Mori's
creation of the theory of extremal rays. This is a crucial principle which
we are only beginning  to exploit \cite {Sa2} \cite{R7} \cite{Co}.

\definition{1.6 Remark}
The minimal model theorem for 3-folds, which, conjecturally, holds
over an arbitrary base scheme, provides a first answer to 1.3, at least if
$\rho (X_K)=1$, namely
a Mori fiber space $X  \to S$ of type 1.4.2.2, that is, one which is a
Del Pezzo fibration. I will now explain why this solution is not satisfactory:

{\bf 1.6.1} To explain this, I need to recall some facts about terminal
singularities of 3-folds (for an excellent introduction to this subject,
see \cite{R5}). If $U$ is a variety with terminal
singularities, the (global) {\it index}
of $U$ is the minimum natural number $I$ such that $IK_U$ is Cartier.
If $p \in U$ is a point, the
index of $U$ at $p$ is the minimum natural number $r$ such that $rK_U$ is
a Cartier divisor near $p$. 3-fold terminal
singularities of index $r=1$
are precisely the isolated $cDV$ singularities (by definition, $p \in U$ is
$cDV$ if the general hyperplane section $p\in B \subset U$ is a surface
with a rational double point at $p$, see 2.3, 2.4)
\cite{R2}. There is a detailed classification of terminal 3-fold singularities
\cite{Mo2}
and the typical example is (up to analytic equivalence):
$$\bigl(p\in U \bigr)= \bigl(0\in (xy+f(z^r,t)=0)\bigr)$$
in affine toric 4-space ${1 \over r}(1, -1, a ,0)$. An isolated singularity
of this form is always terminal of index $r$ if $(a,r)=1$.

The point here is that there are 3-fold Mori fiber spaces over a curve with
singular points of arbitrarily large index. In particular, the minimal
embedding dimension of these models can be arbitrarily large, which is quite
unpleasant, especially for arithmetic purposes. This phenomenon is
very common, even when the generic fiber is ${\Bbb P}^2$ \cite{Mn}.

{\bf 1.6.2} A smooth Del Pezzo surface
$X_K$ of degree $d \geq 3$ is anticanonically
embedded in ${\Bbb P}^d_K$, and similar properties are known for Del Pezzo
surfaces of degree 2 and 1 (3.13).
Returning to our original question, it is natural to look for integral models
$X \subset {\Bbb P}^d_S$ embedded in the same projective space as
$X_K$. This is impossible, in general, for Mori fiber spaces, as
already noted.

{\bf 1.6.3} Fibers of $X\to S$ are often not reduced.
This is also a very common
phenomenon. In particular, not all fibers are Del Pezzo surfaces in the
sense of 1.1. It is very difficult to control the fibers which may
appear on an arbitrary Mori fiber space. It would be nice to have
distinguished integral models whose fibers are completely understood.
\enddefinition

On the other hand, for a variety of reasons, Mori fiber spaces are the
natural objects to look at. Our guiding principle
will be therefore to look for a Mori fiber space with singularities of
{\it smallest possible index}. In most cases, it will be possible to construct
integral models with terminal singularities of index 1. Before I state the
main results of this paper, I recall the following statement due to
Sarkisov \cite{Sa1}, which answers, for conic bundles,
the problem analogous to 1.3:

\proclaim{1.7 Theorem} Let $X \to T$ be a Mori fiber space with
$\operatorname{dim}(X)=\operatorname{dim}(T)+1$, i.e, $X \to T$ is a conic
bundle. There exists a smooth $T^\prime$, a birational morphism
$T^\prime \to T$, a smooth $X^\prime$, birational to $X$, with a structure
of Mori fiber space $X^\prime \to T^\prime$ fitting into a commutative diagram:
$$
\definemorphism{birto}\dashed \tip \notip
\diagram
X^\prime \dto \rbirto & X \dto \\
T^\prime \rto         & T      \\
\enddiagram
$$
The conditions imply that $X^\prime \to T^\prime$ is a {\rm flat} conic
bundle, the anticanonical system $-K_{X^\prime}$ is relatively very ample,
and defines an embedding $X^\prime \hookrightarrow {\Bbb P}^2_{T^\prime}$
(from these facts, it is easy to completely understand the local structure
of $X^\prime \to T^\prime$).
\endproclaim

1.7 is actually rather easy to prove, and Sarkisov \cite{Sa1} uses it as a
starting point to derive global birational rigidity results for $X^\prime$
in the category of Mori fiber spaces.
{}From various points of view, the result is not quite satisfactory: how large
$T^\prime \to T$ needs to be? Is there a natural ``minimal'' choice
(perhaps allowing some mild singularities)?
\medskip

{\bf b. Statement and discussion of the results.}

Gluing together models over the spectra of the local rings of $S$, it is
enough to answer 1.3 when $S=\operatorname{Spec} {\Cal O}$ is the spectrum
of a discrete valuation ring ${\Cal O}$. To facilitate the
statements of the results, we will from now on assume this to be the case.
We leave it to the reader to formulate the results for arbitrary $S$.
We choose a parameter $t \in {\Cal O}$ and denote $k={\Cal O}/t$ the
residue field. For a scheme $Z$ over ${\Cal O}$, $Z_0$ will denote the
special fiber.
We will also assume, unless otherwise indicated, that
{\it ${\Cal O}={\Cal O}_{C,p}$ is the local ring at $p \in C$ of a smooth
complex curve $C$ defined over the field ${\Bbb C}$ of complex numbers.}

Our results are theorems 1.10, 1.15 and 1.18 below.

Let $X_K$ be a smooth Del Pezzo surface of degree $d$ over $K$.
A {\it model} of $X_K$ is a
scheme $X$, defined and flat over ${\Cal O}$, such that $X_\eta=X_K$.

\definition{1.8 Definition} Assume $d\geq 3$.
A {\it standard model} of $X_K$ is a model $X$
of $X_K$ satisfying the following conditions:

1.8.1 $X$ has terminal singularities of index 1.

1.8.2 The central fiber $X_0$ is reduced and irreducible. In particular
$X_0$ is a Gorenstein Del Pezzo surface.

1.8.3 The anticanonical system $-K_X$ is very ample and
defines an embedding $X \subset {\Bbb P}_{\Cal O}^d$. \enddefinition

\definition{1.9 Remark} The most important case of 1.8 is when
$\rho (X_K)=1$. In this case
$X$ has factorial terminal singularities and $X\to S$
is a Mori fiber space. This
condition alone automatically implies that all fibers are reduced and
irreducible \cite {Mo1}.
\enddefinition

\proclaim{1.10 Theorem}
If $d \geq 3$, a standard model for $X_K$ exists.
\endproclaim

In fact, if $d=3$, we prove more, namely we construct standard models with
isolated $cA_n$, $n \leq 5$, $cD_4$, $cD_5$ or $cE_6$ singularities (see
2.23, and 2.4 for the meaning of $cA_n$, $cD_n$, $cE_n$).

\definition{1.11 Remark} In general, standard models are not unique, but
examples are not easy to construct, see 5.8.
\enddefinition

\definition{1.12 Remark} For cubic Del Pezzo surfaces, \S 2 contains a
program which, conjecturally, should construct distinguished integral models
over (nearly) arbitrary discrete valuation rings. Standard models do not always
exist, but every $X_K$ should have a model $X$ which is either a standard
model, or an {\it exceptional model} (2.9).
\enddefinition

If $d=1$ or $2$, models with terminal singularities of index 1 do not
always exist (see 1.18 below). As the anticanonical
class becomes smaller and smaller, it becomes increasingly difficult to
find the birational transformations necessary to construct standard
models, a fact for which we only have a vague explanation.

\definition{1.13 Definition} Assume $d \leq 2$. A model $X$ of $X_K$ is a
{\it standard model} if:

1.13.1 $X$ has terminal singularities.

1.13.2 The central fiber $X_0$ is reduced and irreducible.

1.13.3 If $d=2$, $2K_X$ is Cartier and very ample. If $d=1$, $6K_X$ is
Cartier and very ample.
\enddefinition

\definition{1.14 Remark} As for 1.8, the most important case of 1.13 is when
$\rho (X_K)=1$. In this case
$X$ has ${\Bbb Q}$-factorial terminal singularities and $X\to S$
is a Mori fiber space. This condition alone automatically implies that all
fibers are irreducible, but not, if $I\geq 2$, reduced. \enddefinition

\proclaim{1.15 Theorem}
If $d=2$, a standard model of $X_K$ exists.
\endproclaim

\proclaim{1.16 Conjecture}
If $d=1$, a standard model of $X_K$ exists.
\endproclaim

Presumably, 1.16 can be proved in the same way as 1.15. The author simply lost
patience with the lengthy analysis (in fact, the proof of 1.15 is already
quite long).
Classifying all fibers that may appear on standard models of degree $2$ or
$1$ and index $I>1$ is a problem which requires a finite amount of work.

\definition{1.17 Remark} If $d=2$ (resp. $d=1$), the anticanonical algebra
defines an embedding $X_K \subset {\Bbb P}(2,1,1,1)_K$ (resp.
$X_K \subset {\Bbb P}(3,2,1,1)_K$). The standard models which we construct
(resp. we expect),
are not, in general, embedded in the corresponding weighted projective
space over ${\Cal O}$, see example 4.11. I don't know if 1.15, say, can be
improved allowing models embedded in ${\Bbb P}(2,1,1,1)_{\Cal O}$ only.
\enddefinition

As I said, if $d \leq 2$, models with terminal singularities of index $1$
do not exist in general. In what follows, $I$ denotes the global index of
the canonical class:

\proclaim{1.18 Theorem} If $d=2$, there exist standard models with $I=2$,
which are not birational to any standard model with $I=1$.
Similarly, if $d=1$, there exist standard models with $I=2$, $I=3$ and $I=6$,
which are not birational to any standard model with $I=1$.
\endproclaim

Presumably, there also exist standard models of global index $I$, which are
not birational to a standard model of index $<I$.

The traditional difficulty
of the Fano-Manin-Iskovskikh-Sarkisov method to \linebreak study birational
groups of
Mori fiber spaces is that, a priori, one is forced to examine linear
systems of arbitrarily large degree. In the proof of 1.18, for obvious
reasons, see 5.1, this difficulty is not present.
It is a fact, however, that in all known examples
everything is determined by low degree (often degree 1). One might hope
that there is a general principle governing the experimental evidence.

\medskip

{\bf c. Open questions.}

I will here mention some of the open questions
that arise naturally in connection with this work. I begin with
questions related to birational geometry.

{\bf 1.19} If $X$ is a ${\Bbb Q}$-Fano 3-fold, and $h^0(-K_X)$ is large
(just how large it really needs to be is not clear), Alexeev \cite{A2}
proves that $X$ is birational to a Fano 3-fold
$U$ with Gorenstein canonical singularities, and the general philosophy and
methods of this paper suggest that $U$ is, in turn, birational to a Fano
3-fold with terminal singularities of index 1.

{\bf 1.20} Improve Sarkisov's result 1.7

{\bf 1.21} Study systematically local birational transformations of Del Pezzo
fibrations, improve our local birational rigidity result 1.18.

{\bf 1.22} Let $X$ be a 3-fold, $X \to {\Bbb P}^1$ a Mori fiber space. Give
rigidity results for $X\to {\Bbb P}^1$ in the category of Mori fiber
spaces, similar to \cite{Sa1} for conic bundles.
There are very few known results about these varieties.

It is well known \cite{Ma} \cite{CT} that if the
generic fiber has degree $\geq 5$, $X$ is birational, over ${\Bbb P}^1$, to
${\Bbb P}^1 \times {\Bbb P}^2$.

If $d=4$, $X$ is always birational to a conic bundle over ${\Bbb P}^1$ and one
can obtain nonrationality criteria this way \cite {A1}.

If $d=3$, the only result known to me is due to Bardelli \cite{Ba},
roughly stating that $X$ is almost never rational.

Strangely, nothing seems to be known explicitly about $d=1$
or $d=2$, which should be the rigidest. The methods of \cite{Ba}
require only small modification to cover these cases also.

Here are some questions related to the arithmetic of rational surfaces:

{\bf 1.23} Prove conjecture 2.11, thereby constructing distinguished
integral models of cubic surfaces
over arbitrary discrete valuation rings. Extend this to other Del Pezzo
surfaces and conic bundle surfaces. Study the case
$\operatorname{ch} k=2, 3$.

{\bf 1.24} If $K$ is a local field, $X_K$ a Del Pezzo surface over $K$ and
$X$ a distinguished integral model of $X_K$, for instance, a standard
model, one wishes to calculate $CH_0(X_K)$ in term of data on $X_0$,
together with the knowledge of the lattice $NS(X_K)$ (the N\'eron-Severi group
of $X_K$) as a
discrete Galois module.
The simplest result of this kind was obtained
by Dalawat \cite{Da}, following unpublished ideas of Bloch.
Essentially, he shows that if $X$
is regular, and $NS(X_K)$ is unramified (this assumption is very strong), the
natural specialization map $sp:CH_0(X_K)\to CH_0(X_0)$ is an
isomorphism. Perhaps, $K$ being fixed, there is only a finite number of
possibilities for
$CH_0(X_K)$? The relevance of this problem has to do, of course, with well
known conjectures of Colliot-Th\'el\`ene and Sansuc \cite{CTS}, describing
the Chow group $CH_0(X_K)$ of a rational
surface, over a global field $K$, in term of the local completions.

{\bf 1.25} Let $X$ be a rational surface, with $\rho=1$, defined
over a number field $K$. Manin and
others \cite{Ma2}, \cite{FMT}  propose some very interesting conjectures,
regarding the growth of the function $\varphi (h)$ counting
$K$-rational points of height bounded by $h$. More precisely, if $X$
contains at least one $K$-rational point,
they conjecture that $\varphi (h)$ is asymptotic to $ch$, for some constant
$c$. Compute $c$ in term of (a?) distinguished integral model.

\medskip

{\bf d. Structure of the paper.}

I will now describe the content of each
section of this paper.

\S 2 contains the main ideas, and the reader who is only interested in the
general principles is invited to read this section only.
A program is given, to construct integral
models of cubic surfaces over arbitrary discrete valuation rings. The
program produces a standard model or an exceptional model
(always a standard model if the residue field is algebraically closed).
The main
problem is that we don't know, at the present, that the program terminates
(we conjecture that it does).
The most serious (but not the only) technical obstacle to proving
termination is resolution of singularities. To be honest, we do require, in
addition, that the residue characteristic is $\not = 2, 3$. It is not
impossible to work with residue characteristic $2, 3$, but I did not
pursue this direction. The rest of the paper is devoted to
circumventing, in the geometric case, the problem with termination, and
extending the framework to all Del Pezzo surfaces.

\S 3. In this section, we prove that standard models of Del Pezzo surfaces
of degree $d \geq 3$ exist, i.e., we prove 1.10.
In other words, here and in the rest of the paper we
work with ${\Cal O}={\Cal O}_{C,p}$, the local ring at $p$ of a normal
complex curve $C$. We introduce and study Gorenstein anticanonical models
(many of our results are due, independently, to Alexeev \cite{A2}, but our
proofs are different). The construction of these models uses the theory of
minimal models for 3-folds. If $X_K$ is a Del Pezzo surface of degree $3$, and
$X$ a Gorenstein anticanonical model of $X_K$, it is very easy to see that,
when feeding $X$ into the program of \S 2, termination holds. The $d>3$
case is a fairly easy consequence of the $d=3$ case, together with well
known properties of Del Pezzo surfaces over function fields, and some ideas
from the minimal model program for 3-folds.

In \S 4 we prove 1.15. The ideas are simple variations on the ideas in \S 2
and \S 3, but the actual proof is quite long.

In \S 5 we prove the local rigidity result 1.18, and give an example (5.8)
showing that standard models are not unique.

\medskip

{\bf e. Acknowledgments.} This paper is a considerably revised version of
my 1992 University of Utah doctoral dissertation, written under
the supervision of Prof. J\'anos Koll\'ar, to whom I am deeply indebted,
for too many reasons to list.

In December 1992, I circulated a preprint claiming the results proved here.
In December 1993, S. Mori pointed out an error in a crucial place, and,
a bit later, another error was found.
The present form of the paper benefits from the discovery of a more
conceptual approach to the questions studied, one which allows to considerably
simplify and clarify the proof.

It is largely due to the encouragement of J. Koll\'ar, S. Mori and M. Reid,
that this work now appears, and I wish to express my sincere thanks to all
of them.

Many people helped me with the 1992 version, among them
V. Alexeev, A. Fujiki, Y. Kawamata, J. Koll\'ar, S. Mori and M. Reid.
I am grateful to them all.

When I first submitted the 1992 version, I received valuable comments from the
referees. They are responsible for whatever improvement can be observed in the
presentation of the material.

It is my pleasure to acknowledge a discussion with B. Kunyavskii, who
made me realize some philosophical points expressed in the introduction.

\newpage

\head 2. Integral models of cubic surfaces \endhead

In this section, I will describe a conjectural procedure to
construct distinguished integral models of cubic surfaces over (nearly)
arbitrary discrete valuation rings. This will be used, in \S 3, to actually
construct models over ${\Bbb C}$. After setting up the notation, I will
describe the contents of this section in more detail.

\definition{2.1 Notation}
We fix an arbitrary discrete valuation ring ${\Cal O}$, with
fraction field $K$, parameter $t$, and residue field $k$. My favorites are
the local ring ${\Cal O}_{C,p}$ of a smooth complex  curve $C$ at a point
$p\in C$, ${\Bbb C}\{t\}$, ${\Bbb Q}[t]$, ${\Bbb F}_p[t]$, ${\Bbb
F}_p[[t]]$, the formal
completion ${\Bbb Z}_p$ of the integers at $p$.
We will assume that $\operatorname{ch} k \not =2, 3$.
$\overline k$ is the {\it separable} algebraic closure of $k$, $\overline K
\supset K$ an unramified extension of $K$ inducing $\overline k$,
$\overline{{\Cal O}}$ the integral closure of ${\Cal O}$ in $\overline K$.
$\overline{{\Cal O}}$ is a discrete valuation ring with parameter $t$.
We denote $S = \operatorname{Spec} {\Cal O}$, $\eta=
\operatorname{Spec} K$, $0=  \operatorname{Spec} k$. In what follows all
varieties, schemes, morphisms, etc. are tacitly assumed to be defined over
$S$, unless otherwise indicated. For a scheme $Z$, $Z_\eta$, $Z_0$ are
its generic and special fiber.
$\overline S = \operatorname{Spec} \overline{{\Cal O}}$, and we denote
$\overline Z$ the base change to $\overline S$, $\overline Z_\eta$,
$\overline Z_0$ its generic and special fiber.
A {\it birational map} is always assumed
to be biregular when restricted to the generic fibers.
\enddefinition

The material in this section is organized as follows.
After some preliminaries on $cDV$ and
terminal singularities, we define {\it standard} (2.8) and {\it
exceptional} (2.9) models, and conjecture (2.11) that every cubic Del Pezzo
surface over $K$ has an integral model which is standard or exceptional.
Then (2.13) we describe the flowchart of a program to construct
distinguished integral models, we conjecture (2.14) that the program
terminates, and we state (2.15) that the final product of the program (assuming
that it terminates) is a standard or exceptional model. Conjecture 2.11
then follows from the flowchart 2.13, theorem 2.15, and conjecture 2.14.
In lemmas 2.17--19 we show that each step of the program 2.13 increases
discrepancies, in other words improves singularities. Discrepancies are
rational numbers, so this is not enough to prove 2.14, but it allows to
prove a special case (2.20), which will be instrumental in the proof of
theorem 1.10 (this is done in \S 3). Finally, in the remaining part of this
section, we prove 2.15, using the method of \cite{BW}. The crucial
technical result is 2.26, which permits
to recognize $cDV$ singularities from their equations (see the
``recognition principle'' \cite{BW}, corollary on page 246).

\definition{2.2 Definition} Let $Z_K$ be a scheme over $\operatorname{Spec}
K$. An {\it integral model}, or simply a  {\it model} of $Z_K$ is a
scheme $Z$ over $\operatorname{Spec}
{\Cal O}$ with generic fiber $Z_\eta =Z_K$.
\enddefinition

Given a smooth cubic surface $X_K \subset {\Bbb P}^3_K$, we wish to
construct a nice integral model $X \subset {\Bbb P}^3_{\Cal O}$ of $X_K$.
Before introducing standard (2.8) and exceptional (2.9) models, we recall the
definitions of $cDV$ and terminal singularities.

\definition{2.3 Definition} A normal 2-dimensional scheme $B$ has Du Val
singularities if it has rational singularities of multiplicity 2.
\enddefinition

Over ${\Bbb C}$ (resp. an algebraically closed field of characteristic
$\not = 2, 3, 5$) a surface Du Val singularity is analytically
(resp. formally) equivalent to a standard $A_n$, $D_n$ or $E_6$, $E_7$,
$E_8$ singularity (see for instance the method described in \cite{AVGZ},
vol.1, 12.6, theorem on page 174). In general,
it is still true that the dual graph in the minimal resolution is an
$A_n$, $D_n$ or $E_6$, $E_7$, $E_8$ graph, and we say that the
singularity is of type $A_n$, $D_n$, $E_n$ according
to its graph. References for the general case are \cite{Ar} \cite{Li}.

\definition{2.4 Definition} A 3-dimensional scheme $X$ over $S$ has
$cDV$ (compound Du Val) singularities if, for every singular
$\overline k$-rational point $p
\in \overline X$, there is a hyperplane section $\overline B \ni p$ with a
Du Val singularity at $p$. We say that $p \in X$ is a singularity of type
$cA_n$, $cD_n$ or $cE_n$, according to the type of a general hyperplane
section.
\enddefinition

In the following definitions 2.5, 2.7, a {\it variety} is a reduced and
irreducible scheme over a field, or a reduced and irreducible scheme,
defined and flat over $S$.

\definition{2.5 Definition-Proposition} Let $(X, B)$ be a pair consisting of
a normal variety $X$ and a ${\Bbb Q}$-Weil divisor $B \subset X$ (in other
words, $B=\sum b_iB_i$, where $b_i\in {\Bbb Q}$ and $B_i$ are prime divisors).
Assume that $K_X+B$ is ${\Bbb Q}$-Cartier.

2.5.1 Let $\nu$ be
a discrete valuation of the function field $\operatorname{Frac} X$
of $X$, with center on $X$. By a result of Zariski (see \cite{Ar}),
there exists a projective birational morphism $f:Z \to X$,
from a normal variety
$Z$, such that $\nu=\nu_E$ is the valuation associated to a (prime) divisor
$E\subset Z$. For such an $f$, write:
$$K_Z=f^\ast (K_X+ B)+aE+\sum a_iD_i,$$
where $a_i\in {\Bbb Q}$ and $E$, $D_i$ are distinct prime divisors.
The (rational) number $a$ is independent of $f:Z\to X$, and it is called the
{\it discrepancy} of $\nu$ relative to $K_X+B$ and
denoted $a=a(\nu,K_X+B).$

2.5.2 We say that $X$ has terminal, resp. canonical singularities if
$a(\nu, K_X)>0$, resp. $a(\nu, K_X)\geq 0$, for all valuations $\nu$
with small center on $X$, resp. with center on $X$. \qed
\enddefinition

\definition{2.6 Remark} It is easy to check, by explicit resolution, that
isolated $cDV$ singularities are terminal. In dimension 3 and
characteristic 0, terminal singularities with Cartier canonical class are
isolated $cDV$ singularities \cite{R5}. This is not known in positive or mixed
characteristic, the point being that it is not known, in positive or mixed
characteristic, that terminal singularities are Cohen-Macaulay.
\enddefinition

\definition{2.7 Definition} A variety $X$ is ${\Bbb Q}$-factorial if every
Weil divisor on $X$ is ${\Bbb Q}$-Cartier. Warning: this notion is not
local in the analytic (or \'etale) topology.
\enddefinition

We now define the distinguished integral models which we will try to construct.

\definition{2.8 Definition} Let $X_K \subset {\Bbb P}^3_K$ be a smooth
cubic surface. A subscheme $X \subset {\Bbb P}^3$, flat over $S$, is
said to be a {\it standard integral model} for $X_K$ if $X_\eta = X_K$ and:

2.8.1 $X$ has isolated $cDV$ singularities,

2.8.2 $X_0$ is (reduced and) $k$-irreducible.

In particular, $X$ has Gorenstein terminal singularities and, if
$\operatorname {Pic} X_K={\Bbb Z}$, $X$ necessarily has ${\Bbb
Q}$-factorial singularities and $X \to S$ is a Mori fiber space.
\enddefinition

\definition{2.9 Definition} Let $X_K \subset {\Bbb P}^3_K$ be a smooth
cubic surface. A subscheme $X \subset {\Bbb P}^3$, flat over $S$, is
said to be an {\it exceptional integral model} for $X_K$
if $X_\eta = X_K$ and:

2.9.1 $\overline X_0 =L_1+L_2+L_3$ is union of 3 planes, none of
which is defined over $k$.

2.9.2 $X$ is singular along $C=(L_1\cap L_2) + (L_1\cap L_3)+(L_2\cap
L_3)$.

2.9.3 Let $p \in X_0$ be the triple point. Then $X$ has multiplicity $\mu=2$
at $p$.
\enddefinition

\definition{2.10 Remark} Using the techniques of \S 5, it is possible to
show that an exceptional model is {\it never} birational to a standard
model (we omit the verification).
\enddefinition

Note that exceptional models exist only if $k$ is not algebraically closed.
The conditions imply that $X$ has $cDV$ singularities, $cA_1$ along $C$ and
away from $p$, $cD_4$ at $p$.
A model with terminal singularities can be obtained, from an
exceptional model, by blowing up
$C$ and contracting the strict transforms of the $L_i$s. The resulting
has three index 2 geometric closed points, permuted by the Galois
group $Gal (\overline k/k)$.
\medskip

Our purpose in this section is trying to prove the following:

\proclaim{2.11 Conjecture} Let $X_K$ be a smooth cubic surface over
$K$. Then a model $X$ of $X_K$ exists, which is a standard model or an
exceptional model.
\endproclaim

The main result in this section states that conjecture 2.11 follows from
conjecture 2.14 below. I will soon describe a flowchart for a program to
construct standard or exceptional models. In order to do so, I must first
discuss elementary transformations of projective space.

\medskip

\noindent {\bf 2.12 Elementary transformations of projective space.}
Let ${\Bbb P}={\Bbb
P}^n_{\Cal O}$ be $n$-dimensional projective space over $S$, $L=L_d
\subset {\Bbb P}_0$ be a $d$-dimensional linear subspace, defined over $k$.
If $d \leq n-1$, there is a birational transformation $\Phi=\Phi_L:
{\Bbb P}\dasharrow {\Bbb P}$, centered at $L$. This is nothing but
projection from $L$, and, in homogeneous coordinates,
$\Phi: (x_0:\cdots :x_n)\to
(tx_0: \cdots :tx_d:x_{d_1}:\cdots :x_n)$. $\Phi$ fits into a commutative
diagram:
$$
\definemorphism{birto}\dashed \tip \notip
\diagram
       &A \dlto_f \drto^g &        \\
{\Bbb P} \rrbirto^\Phi&       & {\Bbb P}^+\\
\enddiagram
$$
where $f=bl_L {\Bbb P}: A=Bl_L {\Bbb P}\to {\Bbb P}$ is the blow up of
$L \subset {\Bbb P}$. We denote $F$, $G\subset A$ the $f$ and
$g$-exceptional divisors, so that $f(G)={\Bbb P}_0$ and $g(F)={\Bbb
P}^+_0$. Clearly:
$$\align
K_A &= f^\ast  (K_{\Bbb P}) + (n-d)F \tag 2.12.1 \\
    &=g^\ast   (K_{{\Bbb P}^+})+(d+1) G \tag 2.12.2
\endalign$$
I will often simply denote $K, K^+$ the canonical classes of ${\Bbb P},
{\Bbb P}^+$.

One can show that any birational $\Phi : {\Bbb P}_{\Cal O}
\dasharrow {\Bbb P}_{\Cal O}$ is a composition of $\Phi_L$s, with $L
\subset {\Bbb P}_0$ of dimension $n-1$, and a biregular isomorphism.

\medskip

We are now ready to describe the (still conjectural) procedure to construct
distinguished integral models of cubic surfaces. Fix a smooth cubic
surface $X_K \subset {\Bbb P}^3_K$.
\medskip
\noindent {\bf 2.13 Flowchart.}

{\bf 2.13.0} Let $X \subset {\Bbb P}={\Bbb P}^3_{\Cal O}$ be an arbitrary
flat closure of $X_K$.

{\bf 2.13.1} Does $X$ have generic multiplicity $\mu \geq 2$ along a 2-plane
$L\subset {\Bbb P}_0$? If yes, let $\Phi_L$ be as in 2.12, $X^+=\Phi_{L
\ast}X$, and go back to 2.13.1 with $X^+$ in place of $X$. If not, go to
2.13.2.

{\bf 2.13.2} Does $X$ have generic multiplicity $\mu=3$ along a line
$L\subset {\Bbb P}_0$? If yes, let $\Phi_L$ be as in 2.12, $X^+=\Phi_{L
\ast}X$, and go back to 2.13.1 with $X^+$ in place of $X$. If not, go to
2.13.3.

{\bf 2.13.3} Does $X_0$ contain a plane $L$? If yes, let $\Phi_L$ be as in
2.12,
$X^+=\Phi_{L\ast}X$, and go back to 2.13.1 with $X^+$ in place of $X$. If
not, go to 2.13.4.

{\bf 2.13.4} Does $X$ have generic multiplicity $\mu=2$ along a line
$L\subset {\Bbb P}_0$? If yes, let $\Phi_L$ be as in 2.12, $X^+=\Phi_{L
\ast}X$, and go back to 2.13.1 with $X^+$ in place of $X$. If not, go to
2.13.5.

{\bf 2.13.5} Is there a $k$-rational point $p \in X$ of multiplicity 3?
If yes, let $\Phi_L$ be as in 2.12, $X^+=\Phi_{L
\ast}X$, and go back to 2.13.1 with $X^+$ in place of $X$. If not,
by theorem 2.15 below, one of
the following is true:

{\bf 2.13.5.1} $X$ is a standard model.

{\bf 2.13.5.2} $X$ is an exceptional model.
\medskip

By the flowchart just described, conjecture 2.11 follows from the
following conjecture 2.14, and theorem 2.15:

\proclaim{2.14 Conjecture} The program described by the flowchart terminates.
\endproclaim

\proclaim{2.15 Theorem} Let $X\subset {\Bbb P}^3_{\Cal O}$ be a subscheme,
flat over ${\Cal O}$, whose generic fiber $X_\eta \subset {\Bbb P}^3_\eta$
is a smooth cubic surface. Assume the following:

2.15.1 $X_0$ is $k$-irreducible. This is equivalent to saying that $X_0$
contains no 2-plane defined over $k$, and it implies that $X_0$ is reduced.
In particular, $X$ is nonsingular in codimension one, hence normal.

2.15.2 $X$ is nonsingular at the generic point of every line $L
\subset {\Bbb P}_0$ defined over $k$.

2.15.3 $X$ has multiplicity $\mu \leq 2$ at every $k$-rational point $p \in
X$.

Then either $X$ is a standard model, or $X$ is an exceptional model.
\endproclaim

The rest of this section is devoted to two tasks.
First, I will explain the philosophy behind 2.14 and prove a special case of
termination(2.20). This special case,
together with some ideas from the theory of
minimal models, will be used, in \S 3, to prove the existence of
standard models in the geometric setting. Finally, I will prove 2.15.

\medskip

I now begin discussing conjecture 2.14. It is customary, in the theory of
minimal models of algebraic varieties, to measure singularities of a
variety $U$ according to the discrepancies $a(\nu, K_U)$, taken with
respect to the canonical class $K_U$, of the valuations $\nu$ with small
center on $U$: the higher the discrepancies, the better the singularities.
Lemmas 2.17--19 below state that each transformation $\varphi_L: X \dasharrow
X^+$ used in 2.13 increases discrepancies, in other words passing from $X$
to $X^+$ {\it improves} singularities. This implies an important special
case of termination 2.20. In general, one is tempted to speculate that any
process that consistently increases discrepancies must eventually stop.
The following is elementary:

\proclaim{2.16 Lemma-definition}
Let $U$ be a variety with canonical singularities.
A valuation $\nu$, of the fraction
field $\operatorname{Frac} U$ of $U$, with small center
on $U$, is {\rm crepant} if $a(\nu, K_U)=0$. Assume there is a nonsingular
variety $Z$, and a proper birational morphism
$f:Z \to U$. Then, the number of crepant valuations of $U$ is finite,
and is denoted $e(U)$. \qed
\endproclaim

In the following statements 2.17--20, we fix a smooth cubic surface $X_K
\subset {\Bbb P}^3_K$ and a model $X \subset {\Bbb P}_{\Cal O}$ of $X_K$.
$n(X)$ denotes the number of $k$-irreducible components of the central
fiber $X_0$.

\proclaim{2.17 Lemma} (Blowing up the plane) assume $X_0$ contains a 2-plane
$L\subset {\Bbb P}_0$, defined over $k$. Let $1 \leq \mu \leq 3$ be the
generic multiplicity of $X$ along $L$.  Let $\Phi=\Phi_L: {\Bbb
P}\dasharrow {\Bbb P}^+$, $X^+=\Phi_\ast X$ be as in 2.13:

2.17.1 Let $\nu$ be a valuation of $K({\Bbb P})$. Then:
$$a\bigl(\nu, K^++X^+\bigr)\geq a\bigl(\nu, K+X+(\mu-1){\Bbb P}_0\bigr).$$

2.17.2 If $X$ has canonical singularities, so does $X^+$, and
(assuming resolution of singularities) $e(X^+)<e(X)$, or $e(X^+)=e(X)$ and
$n(X^+)<n(X)$.
\endproclaim

\proclaim{2.18 Lemma} (Blowing up the line) assume $X$ has generic
multiplicity
$2 \leq \mu \leq 3$ along a line
$L\subset {\Bbb P}_0$, defined over $k$. Let $\Phi=\Phi_L: {\Bbb
P}\dasharrow {\Bbb P}^+$, $X^+=\Phi_\ast X$ be as in 2.13:

2.18.1 For $\nu$ as in 2.17.1:
$$a\bigl(\nu, K^++X^+\bigr)\geq a\bigl(\nu, K+X+(\mu-2){\Bbb P}_0\bigr).$$

2.18.2 Same as 2.17.2.
\endproclaim

\proclaim{2.19 Lemma} (Blowing up the point) let $p\in X$ be a $k$-rational
point of multiplicity $\mu =3$. Let $\Phi=\Phi_{\{p\}}: {\Bbb
P}\dasharrow {\Bbb P}^+$, $X^+=\Phi_\ast X$ be as in 2.13:

2.19.1 For $\nu$ as in 2.17.1:
$$a\bigl(\nu, K^++X^+\bigr)\geq a\bigl(\nu, K+X\bigr).$$

2.19.2 Same as 2.17.2.
\endproclaim

\proclaim{2.20 Corollary} (Assuming resolution of singularities) if $X$ has
canonical singularities, starting the flowchart 2.13 with $X$ in 2.13.0,
termination holds.
\endproclaim
\demo{Proof} Lemmas 2.17--19 and descending induction on $e(X), n(X)$.\qed
\enddemo

The above results will follow easily from the following result. The
formulation is a bit technical, but it is useful, and will be used in this
form in the construction of integral models of Del Pezzo surfaces of degree
2, in \S 4.

\proclaim{2.21 Lemma} Let $p:U \to S$, $p^+:U^+ \to S$ be proper
morphisms. Assume that $K_U$, $K_{U^+}$ are ${\Bbb Q}$-Cartier, and $-K_U$
is $p$-ample. Let $\Phi : U \dasharrow U^+$ be a birational map defined
over $S$ (i.e, $p=p^+ \Phi$). Assume that $a(\nu_F, K_U) \leq 0$ for every
$\Phi^{-1}$-exceptional divisor $F \subset U^+$. Then:

2.21.1 $a(\nu, K_U) \leq a(\nu, K_{U^+})$ for all valuations $\nu$
(exceptional or
not). In particular, if $U$ has canonical singularities, so does $U^+$ and,
assuming resolution of singularities, $e(U^+) \leq e(U)$.

2.21.2 $a(\nu_G, K_{U^+})>0$ for every $\Phi$-exceptional
divisor $G \subset U$.

2.21.3 Assuming resolution of singularities, if $U$ has canonical
singularities, $e(U)=e(U^+)$, and there are no
$\Phi$-exceptional divisors $G\subset U$, $\Phi$ is an isomorphism in
codimension 1. If, in addition, $-K_{U^+}$ is $p^+$-ample, $\Phi$ is
an isomorphism.
\endproclaim

\demo{Proof} We prove 2.21.1 first. Let $\nu$ be a valuation, and
$\Gamma$ the normalization of the graph of $\Phi$. By a theorem of
Zariski (see \cite{Ar}),
there is a normal variety $Z$ and a projective birational morphism
$Z \to \Gamma$ such that $\nu$ is a divisor in $Z$.
Let $f: Z \to U$, $g: Z \to U^+$ be the natural maps. We will prove that
$a(\nu_E ,K_U) \leq a(\nu_E , K_{U^+})$ for every divisor $E \subset Z$.
This will clearly prove 2.21.1.
Let $E_i \subset Z$ be the divisors which are exceptional over both $U$ and
$U^+$, $G_j$ the divisors which are exceptional for $g$ but not for $f$,
and similarly $F_k$ the divisors which are exceptional for $f$ but not for
$g$. We have:
$$\align
K_Z &= g^\ast K^+ + \sum a_i^+E_i+\sum b_j^+ G_j \\
    &= f^\ast K   + \sum a_i  E_i+\sum c_k   F_k.
\endalign$$
Then, denoting $\equiv_g$ numerical equivalence with respect to $g$, we
have:
$$-\sum b_j^+G_j -\sum (a_i^+-a_i)E_i \equiv_g -f^\ast K - \sum c_kF_k, \eqno
(\ast)$$
where by assumption $-f^\ast K$ is $g$-nef and $-c_k \geq 0$.
By the negativity lemma (see e.g. \cite{Co}, 2.5), $ b_j^+ \geq 0$ and $a_i^+
\geq a_i$ for all $j$, $i$. This concludes 2.21.1.

We now prove that all $b_j^+>0$. Assume by contradiction that $b_1^+=0$,
say. Restricting $(\ast )$ to $G_1$ we get:
$$-\sum b_j^+G_j|_{G_1} -\sum (a_i^+-a_i)E_i|_{G_1} \equiv_g -
f^\ast K|_{G_1} - \sum c_kF_k|_{G_1}.$$
The divisor on the left is antieffective, the divisor on the right is the
sum of a nef and big divisor and an effective divisor. This is a
contradiction, which proves 2.21.2.

Assume now that $U$ has canonical singularities (hence $U^+$ does too), and
let $Z$ as above be a resolution of $\Gamma$.
By assumption all $c_k=0$, so the divisors contributing to
$e(U)$ are the $F_k$s and the $E_i$s with $a_i=0$.
By 2.21.2, the divisors contributing to $e(U^+)$ are the $E_i$s with
$a_i^+=0$, and, since $a_i^+\geq a_i\geq 0$, $e(U)=e(U^+)$ implies there are no
$F_k$s. Therefore, if there are no $G_j$s, $\Phi$ is an isomorphism in
codimension 1. The last part of 2.21.3 is also easy, and is a special case
e.g. of \cite{Co}, 2.7. \qed \enddemo

The proofs of 2.17, 2.18, 2.19 are very similar and easy, so I will only
prove 2.17.

\demo{Proof of 2.17} Let $A$, $f: A \to {\Bbb P}$, $g: A \to {\Bbb P}^+$,
$F$ and $G$ be
as in 2.12, $Z=f^{-1}_\ast X \subset A$. First of all let us prove
that: $$g^\ast X^+=Z+(3-\mu)G.$$ Indeed, in a neighborhood of the
generic point of $L$, $X$ is
given by an equation $x^k+t\varphi (t, x)=0$. The relevant chart for the
blow up is $t=xt^\prime$, so $x^k+t\varphi (t, x)=x^k+xt^\prime \varphi
(xt^\prime, x)= x^\mu\bigl(x^{k-\mu}+t^\prime \varphi^\prime (t^\prime, x)
\bigr)$. Then $g^\ast X^+=Z+\bigl((3-k)+(k-\mu)\bigr)G=Z+(3-\mu)G$.

The crucial formulas are:
$$\align
f^\ast\bigl(K+X\bigr)&=K_A+Z+(\mu-1)F, \\
g^\ast\bigl( K^++X^+\bigr)             &=K_A+Z-\mu G.
\endalign$$
This shows that the assumptions of 2.21 are satisfied with $X=U$,
$X^+=U^+$, so the first statement in 2.17.2 follows from 2.21.1. It is easy
to see that $X \dasharrow X^+$ is not an isomorphism in codimension one: if
$\mu \geq 2$, there is a divisor exceptional for $Z \to X$ but not for $
Z \to X^+$, if $\mu=1$, there is a divisor exceptional for $Z \to X^+$ but
not for $Z \to X$. The rest of 2.17.2 then follows from 2.21.2 and
2.21.3. To prove 2.17.1, observe that our crucial formulas imply that
$a({\Bbb P}^+_0, K+X+(\mu -1){\Bbb P}_0 )\leq 0$. 2.17.1 can then be proved
with the same technique as 2.21.1 (left to the reader). \qed \enddemo

This concludes the discussion of conjecture 2.14. The rest of this section
is devoted to the proof of theorem 2.15. We must show that the end result
$X$ of the flowchart 2.13 is an exceptional or standard model.
We first show that either $X$ is an exceptional model, or it has isolated
singularities. As a second step, we show that $X$ is a standard model if it
has isolated singularities:

\proclaim{2.22 Lemma} Let $X$ be as in 2.15. If the singular set of $X$
contains a curve, $X$ is an exceptional model.
\endproclaim

\demo{Proof} $X_0$ is reduced and $k$-irreducible, so $\overline X_0$
is reduced. If $X_0$ is not geometrically
irreducible, $\overline X_0=L_1+L_2+L_3$ is the union of three 2-planes
$L_i$, none
of which is defined over $k$: indeed, if $\overline X_0=L+Q$, where $L$ is
a 2-plane and $Q$ an irreducible quadric, the 2-plane $L$ is necessarily
defined over $k$.

If $X_0$ is geometrically irreducible, the nonnormal locus of $X_0$, if
nonempty, consists of a line. Then $X$ has isolated singularities by
condition 2.15.2.

Otherwise, $\overline X_0=L_1+L_2+L_3$ as above. Let
$C=(L_1\cap L_2) + (L_1\cap L_3)+(L_2\cap L_3)$. No component of $C$ is
defined over $k$, for otherwise one of the $L_i$s is defined over $k$:
for instance if $L_1 \cap L_2$ is defined over $k$, $L_3$ is also
necessarily defined over $k$. This means that either $X$ has isolated
singularities, or the singular locus of $X$ is all of $C$. But in this
case, since $p=L_1 \cap L_2 \cap L_3$ is necessarily $k$-rational,
$X$ is an exceptional model.
\qed \enddemo

2.15 is an immediate consequence of 2.22 and 2.23 below:

\proclaim{2.23 Theorem} Let $X \subset {\Bbb P}_{\Cal O}$ be a model of
$X_K$. Assume:

2.23.1 $X_0$ is $k$-irreducible.

2.23.2 $X$ has isolated singularities.

2.23.3 $X$ has multiplicity $\mu \leq 2$ at every $k$-rational point $p \in
X$.

Then $X$ has $cA_n$, $n \leq 5$, $cD_4$, $cD_5$ or $cE_6$ singularities,
in particular $X$ is a standard model.
\endproclaim

\definition{2.24 Remark} 2.23 seems somewhat miraculous to me. For instance,
the conclusion is false if $X_0$ is reducible, and the following is,
essentially, the only
counterexample. Let $X$ be defined by:
$$x_0(x_0x_3+x_1^2)+tx_2^3+t^nx_3^3=0$$
in ${\Bbb P}_{\Cal O}$. $X_0$ is the projective cone over the union of a
smooth conic and a line tangent to it. It is easy to see that $X$ has
isolated singularities, in fact the only singular point is $p=(0,0,0,1)$.
However, if $n \geq 4$, $p\in X$ is not a $cDV$ singularity (it is always a
canonical singularity).
\enddefinition

We will prove 2.23 by explicitly studying the singularities on $X$,
following the method of \cite{BW}.
In characteristic zero, using vanishing theorems, J. Koll\'ar suggested
the following nice conceptual proof that $X$ has
canonical singularities. Isolated
canonical singularities are not necessarily $cDV$, so this is still short
of proving 2.23, although it does provide some conceptual understanding.

\proclaim{2.25 Theorem} (Characteristic zero only) let $X \subset
{\Bbb P}^n_{\Cal O}$ be a hypersurface, of
degree $d \leq n$ and with isolated singularities. Then $X$ has canonical
singularities.
\endproclaim

\demo{Proof} Let $x \in X_0$ be a singular point, $H \subset {\Bbb P}$ be a
general hyperplane, not passing through $x$. Assume $x \in X$ is not
a canonical singularity, choose $\varepsilon, \delta>0$ very small, and
consider the divisor:
$$K_{\Bbb P}+(1-\varepsilon) (X + {\Bbb P}_0)+(1+\delta)H.$$
This divisor is not log canonical at $x$ and not log canonical along $H$,
but, since $X$ has isolated singularities, it is klt in a small
neighborhood of $x$, away from $x$ and $H$. This
contradicts the connectedness of the nonklt locus \cite{FA}, 17.4
(whose proof uses resolution of singularities and vanishing theorems). \qed
\enddemo

We will now prove 2.23. We will be concerned with subvarieties of
$A=$\linebreak $\operatorname{Spec} {\Cal O}[x_1,\dots x_n]$,
defined by an equation:
$$f(x_1, ... ,x_n)=0,$$
where $f \in {\Cal O}[x_1,..., x_n]$ is a polynomial
of degree $\leq 3$ in $x_1, ..., x_n$, with coefficients in ${\Cal O}$.
Abusing language slightly, we will say that the monomial $t^mx^l$ {\it
appears} in $f$, and write $t^mx^l \in f$, if the monomial $x^l$ appears
in $f$ with coefficient $a_l \in {\Cal O}$, and $\operatorname{ord} a_l=m$.
We then write $a_l=u_lt^m$ with $u_l \in {\Cal O}$ a unit, and:
$$f=\sum f_m (x_1, ...,x_n)t^m,$$
where $f_m \in {\Cal O}[x_1, ..., x_n]$ is a polynomial all of whose
coefficients are units in ${\Cal O}$. The {\it order} of $f$ at the origin
is the minimum integer $m+l$ such that $t^mx^l \in f$. If $d$ is the order
of $f$, we write $f_{(d)}$ the part of order exactly $d$, it is the sum of
all monomials $t^mx^l\in f$ with $m+l=d$. Let $f:B \to A$ be the blow up of
the maximal ideal at the origin. The {\it tangent cone} to $f=0$ (at the
origin) is the fiber of $f$ at 0, and the procedure to find it is as follows:
take $f_{(d)}$, replace $t$ with a formal variable $\tau$, reduce mod $t$.
This gives a polynomial $\tilde f_{(d)}(\tau, x_1, ..., x_n)\in
k[\tau, x_1, ..., x_n]$, homogeneous of degree $d$. The tangent cone is
then $(\tilde f_{(d)}=0) \subset \operatorname{Proj}k[\tau, x_1, ..., x_n]$.

The proof of 2.23 is based on the following result, analogous to the
``recognition principle'' \cite{BW}, corollary on page 246, also used
(implicitly) in \cite{R2}, which is proved resolving by explicit blow ups:

\proclaim{2.26 Lemma} Let $B=(f(x,y)=0) \subset A=\operatorname{Spec} {\Cal
O}[x,y]$ be a surface as above. Assume that $B$ has, at the origin, an
isolated singularity of multiplicity 2.

2.26.1 If $0 \in B$ has reduced tangent cone, $0\in B$ is a rational double
point of type $A_n$, for some $n$.

2.26.2 With the above notation, assume that $f_0=f_{(2)}=x^2$ (in
particular, the tangent cone is not reduced). Let
$g(y)=f(0, y)$, assume that $\operatorname{ord} g=3$, and the tangent cone
of $g$ is not a triple point. Then $0\in B$ is a
rational double point of type $D_n$, for some $n$.

2.24.3 Assume that $f_0=f_{(2)}=x^2$ and, with $g$ as in 2.26.2,
$g_0=g_{(3)}=y^3$. If $t^4 \in f$, $0\in B$ is a rational double point of type
$E_6$. \qed \endproclaim

\demo{2.27 Proof of 2.23}
We will use the method of \cite{BW}. First of all,
it is proven in \cite{BW} that a normal cubic surface, which is not a cone,
has rational double points of type $A_n$, $n \leq 5$, $D_4$, $D_5$ or $E_6$
(their proof is easily extended to the case $\operatorname{ch}(k)\not =
2,3$). 2.23 follows if $X_0$ is normal and not a
cone, and, from now on, we assume this not to be the case.

Let us  base change to $\overline {\Cal O}$ to achieve $k$ separably
closed. After base change, $X_0$ is not necessarily irreducible, but it is
reduced and $X_0 \not = L+Q$ with $L$ a plane and $Q$ an irreducible
quadric (recall 2.24), for otherwise $L$ would be Galois invariant and $X_0$
would be reducible before base change.

Let $p \in X$ be a singular point.
Let $p\in B$ be a
hyperplane section of $X$, general among those passing through $p$.
We will prove that $p\in B$ is a Du Val singularity. Let
$ B_0=B\cdot X_0$,
then $ B_0$ is a reduced plane cubic, and $p \in B_0$ is
a singular point. We distinguish five cases:
\medskip
2.27.1 $B_0$ is a nodal rational curve and $p \in B_0$ the node.

2.27.2 $B_0$ is a cuspidal rational curve and $p \in
B_0$ the cusp.

2.27.3 $B_0=\ell+q$ where $\ell$ is a line and $q$ a reduced conic.

2.27.4 $B_0=\ell_1+\ell_2+\ell_3$ and $p\in B_0$ is a
double point.

2.27.5 $B_0=\ell_1+\ell_2+\ell_3$ and
$p \in B_0$ is the triple point.
\medskip
We choose homogeneous coordinates
$x_0,x_1, x_2, x_3$ on ${\Bbb P}$ so that $p=(0, 0, 0, 1)$, then work with
inhomogeneous coordinates $x=x_0/x_3$, $y=x_1/x_3$, $z=x_2/x_3$.
We may assume $B=(z=0)$. If $X_0$ is not normal, since
$X_0$ is not the union of a plane and an irreducible quadric, there always
is a line $r \ni p$ along which $X_0$ is
singular, and we may assume $r=(x=y=0)$. The assumption
$\operatorname{ch} k \not = 2, 3$ will be used to put $B_0$ in a suitable
normal form in each case.

2.27.1 Note that 2.26.1 already implies that $p \in B$ is $A_n$ for some $n$.
We may assume that $B_0$ is described by the equation:
$$x^2+y^2+y^3=0.$$
In this case $X_0$ is not normal and $X_0$, $X$ are described,
respectively, by:
$$\align
&x^2+y^2+y^3+z\varphi(x, y)=0,\\
&x^2+y^2+y^3+z\varphi(x, y)+\sum_{m \geq 1}f_m(x, y, z)t^m=0,
\endalign$$
with $\varphi=\varphi_{(2)}$. Since $X$ is not singular along $r$, $z^l \in
f_1$ for some $l$, and necessarily $l \leq 3$. It is then easy to see that
$n \leq 3$.

2.27.2 We may assume that $B_0$ is described by the equation:
$$x^2+y^3=0.$$
In this case $X_0$ is not normal and $X_0$, $X$ are described,
respectively, by:
$$\align
&x^2+y^3+z\varphi(x, y)=0,\\
&x^2+y^3+z\varphi(x, y)+\sum_{m \geq 1}f_m(x, y, z)t^m=0,
\endalign$$
with $\varphi=\varphi_{(2)}$. Since $X$ is not singular along $r$, $z^l \in
f_1$ for some $l$, and necessarily $l \leq 3$. 2.26 implies that $p \in B$
is $A_{n_1}$, $D_{n_2}$ or $E_6$. A simple analysis shows $n_1 \leq 1$, $n_2
\leq 5$.

2.27.3 This case does not occur: necessarily $X_0=L+Q$ is union of a plane
and a reduced quadric.

2.27.4 This is almost identical to 2.27.1.

2.27.5 In this case $X_0$ is a cone with vertex at $p$. One can easily see
that $p \in X$ is $cA_n$ or $cD_4$ (similar to 2.26.2), and a simple
analysis shows that $n \leq 3$.
\qed \enddemo

\newpage

\head 3. Geometric integral models of Del Pezzo surfaces of degree
$d\geq 3$ \endhead

In this section we prove theorem 1.10. The main point is to only use the
known case 2.20 of termination of the flowchart 2.13.
In order to do this we will construct a model $U$ with Gorenstein canonical
singularities (3.4), and use its relative anticanonical algebra to embed it
in ${\Bbb P}^3$ (3.13). These models, and their properties, are best
understood within the context of ``good models with $-K$ nef'' (3.1).
A different proof of 3.4 is due to Alexeev \cite{A2}.
3.4, 3.13 and 2.20 already imply the existence of standard models of cubic
surfaces. The proof of 1.10 is finally done in 3.18, where the general case is
reduced to the case of cubics, using the cone theorem for models with $-K$
nef (3.16).

The construction of $U$, and the reduction $d=3 \Rightarrow d\geq3$,
utilize ideas from the theory of minimal models of 3-folds, for which
the standard reference is \cite{KMM}. We also use results from
\cite{Ka} and \cite{FA}, but, in an effort to make the presentation
more accessible, we try to clearly state all we need. I believe that good
models with $-K$ nef, together with the cone theorem 3.16, are of
independent interest in the study of models of Del Pezzo surfaces, and may
be useful in other context as well.

We work with the same setup as in \S 2, only now ${\Cal O}={\Cal O}_{C,p}$
is the local ring at $p\in C$ of a smooth complex curve $C$. As in \S 2,
unless otherwise indicated, all varieties, schemes, morphisms, etc., are
assumed to be defined over $S= \operatorname{Spec} {\Cal O}$. If $X$
is a scheme, $X_\eta$, $X_0$ denote the generic and central fiber.
Birational maps are assumed to induce isomorphisms of generic fibers.

\definition{3.1 Definition}
Let $X$ be a scheme over $S$, such that $X_\eta$ is a Del Pezzo
surface. $X$ is a {\it good model with $-K_X$ nef} if:

3.1.1 $X$ has terminal singularities of index 1,

3.1.2 $-K_X$ is nef,

3.1.3 a general member $B \in |-K_X|$ has at worst Du Val singularities
(the conditions imply that $B$ is actually smooth).
\enddefinition

\definition{3.2 Remark}
Let $X$ be a model. Assume that $X$ has
terminal singularities of index 1, and $-K_X$ is nef. I conjecture
that under these assumptions a general member $B \in |-K_X|$ has at
worst Du Val singularities. This should be proven with the
techniques of Kawamata (for an introduction, with results in the context
most relevant to us, I recommend \cite{R4}) and Shokurov \cite{Sh}.
For example, \cite{R4} implies immediately that
$B$ has at worst semi log canonical singularities.
\enddefinition

The following is our main result in this section, and it is due
independently to Alexeev \cite{A2}:

\proclaim{3.3 Theorem}
Let $X_K$ be a Del Pezzo surface. Then, $X_K$ has a
good model $X$ with $-K$ nef.
\endproclaim

\proclaim{3.4 Corollary-Definition}
Let $X$ be a good model with $-K_X$ nef. Then $-K_X$ is eventually free and
defines a morphism $\varphi : X \rightarrow {U} =
\operatorname{Proj} \oplus_{n \geq 0}H^0(-nK_X)$. $U$ has Gorenstein
canonical singularities and $-K_{U}$ is ample.
$U$ is called a {\it Gorenstein anticanonical model}.
Also, ${U}$ is called the {\it anticanonical model}
of $X$.
\endproclaim

\demo{Proof} This is an immediate consequence of the base point free
theorem \cite{KMM}, 3-1-1. \qed \enddemo

It is important to realize that, despite the name, the anticanonical model
is not ``canonical'', i.e. it is not a birational invariant of the original
Del Pezzo surface $X_K$: it depends on the choice (up to flops) of a good model
with $-K$ nef. Later, in 3.13, the structure of the algebra
$\oplus_{n \geq 0} H^0(-nK_X)$ will be studied in detail.

The Gorenstein anticanonical model $U$ determines, up to flop, the
model with $-K$ nef, by a procedure called ``maximal crepant blow up''.
Before I make this precise, let me recall the theory of maximal crepant
blow ups as developed by Kawamata \cite{Ka}.

\proclaim{3.5 Theorem-Definition}
Let $V$ be a 3-fold with canonical singularities.

3.5.1 The following two conditions are equivalent for a birational projective
morphism $\varphi :Z \to V$ from a normal variety $Z$:

3.5.1.1 $K_Z= \varphi^\ast K_V$ and $Z$ has ${\Bbb Q}$-factorial terminal
singularities.

3.5.1.2 $K_Z= \varphi^\ast K_V$ and $Z$ is maximal with this property.

3.5.2 Given $V$, a projective $\varphi : Z \to V$ satisfying the above
equivalent conditions always exists. It is called a
{\rm maximal crepant blow up}
of $V$. Any two maximal crepant blow ups of $V$ differ by a chain
of flops over $V$. In addition, the following properties hold:

3.5.2.1 Every valuation $\nu$ of ${\Bbb C}(V)$,
with discrepancy $a(\nu, K_V)= 0$, is a $\varphi$-exceptional divisor.
In particular
$\rho(Z)-\rho(V) \geq e(V)$, where $e(V)$ is the number of crepant valuations.

3.5.2.2 $\rho(Z)-\rho (V)=e(V)$ if and only if $V$ is
${\Bbb Q}$-factorial.\qed
\endproclaim

\definition{3.6 Remark}
Let $X$ be a good model with $-K$ nef, $U$ its anticanonical
model. Then $\varphi : X \to U$ is a maximal crepant blow up. It is easy to
see that every other maximal crepant blow up $\varphi^\prime : X^\prime \to
U$ is a good model with $-K_{X^\prime}$ nef.
\enddefinition

The first ingredient in the proof of 3.3 is the following result:

\proclaim{3.7 Theorem}

3.7.1 Let $K$ be a $C_1$ field (the field $K=k(C)$ of rational functions of an
algebraic curve $C$ defined over an algebraically closed field is an
example of a $C_1$ field), and $X$ a Del Pezzo surface over $K$. Then,
the set $X(K)$ of $K$-rational points of $X$ is not empty.

3.7.2 Let $K=k(C)$ be the field of rational functions of an
algebraic curve $C$ defined over an algebraically closed field. Then,
the set $X(K)$ of $K$-rational points of $X$ is
dense in the Zariski topology of $X$.
\endproclaim

\demo{Proof}
3.7.1 is in Colliot-Th\'{e}l\`{e}ne
\cite{CT}. It is usually possible to prove directly, using
geometric arguments, that the $K$-rational points are dense in
the Zariski topology. 3.7.2 is \cite{KoMM}, 2.13. \qed
\enddemo

Note that, over any field $K$, all Del Pezzo surfaces $X$ of degree 1 have
one distinguished $K$-rational point, the base locus of the anticanonical
linear system $|-K_X|$.

The second ingredient in the proof of 3.3 is a special case of the
logarithmic minimal model program for 3-folds, which I now recall
following \cite{FA}, Ch 5. We only outline (some of)
the proofs. For a complete treatment, the reader is referred to loc.cit.

\definition{3.8 Definition}
Let $X$ be a 3-fold and $B \subset X$ a reduced divisor.
The pair $(X, B)$ satisfies condition $(\ast)$ if:

3.8.1 $X$ is normal and ${\Bbb Q}$-factorial,

3.8.2 for all normal varieties $Z$ and birational morphism $f: Z \to X$
with exceptional divisors $E_i$ one has:
$$K_Z + f^{-1}_\ast B=f^\ast (K_X+B)+\sum a_i E_i$$
with all $a_i \geq 0$ and $a_i=0$ only if the center of $E_i$ on $X$
is a curve contained in $B$.
\enddefinition

The following is the crucial consequence of condition $(\ast)$:

\proclaim{3.9 Lemma}
Assume $(X, B)$ satisfies condition $(\ast)$. Then:

3.9.1 $X$ has terminal singularities.

3.9.2 If $x \in B$, $X$ and $B$ are smooth at $x$.
\endproclaim

\demo{Proof}
The first statement is obvious. Let $S \subset B$ be a component containing
$x$. Then $S, \operatorname{Diff}_S(B-S)$ is terminal, so $S$ is smooth and
$\operatorname{Diff}=0$. Let $\pi:
x^\prime, X^\prime \to x, X$ be the index 1 cover in a neighborhood of $x$.
This covering is \'etale outside $x$, it has degree $d$, the index of
$K_X$ at $x$, and it is totally ramified at $x$. $S$ is smooth so
$\pi^{-1}S$ is the union of $d$ irreducible components meeting at
$x^\prime$. $\pi^{-1}(S)$ is a ${\Bbb Q}$-Cartier divisor on a $cDV$
3-fold, so it is Cartier. Therefore $\pi^{-1}S$ satisfies Serre's
condition $S_2$, and letting $i:U =
\pi^{-1}(S)\setminus \{x^\prime\}\to \pi^{-1}(S)$ be the injection,
$i_\ast {\Cal
O}_U= {\Cal O}_{\pi^{-1}S}$. This implies that $U$ is connected, hence
$\pi^{-1}S$ is irreducible, so $d=1$, $x \in X$ is a $cDV$ point and
$S$ is Cartier. Since $S$ is smooth, $X$ is also smooth.
\qed \enddemo

\proclaim{3.10 Theorem} Let $(X, B)$ satisfy condition $(\ast)$.
Let $f: X \to Z$ be a small extremal contraction such that
$-(K_X+B)$ is $f$-ample. Then, the flip of $f$ exists. \qed
\endproclaim

\proclaim{3.11 Proposition} Let $(X, B)$ satisfy condition
$(\ast)$. Let $(X^\prime, B^\prime)$ be obtained from $(X, B)$ by a
sequence of $K+B$-flips and divisorial contractions not contracting a
component of $B$. Then $(X^\prime,B^\prime)$ satisfies condition $(\ast)$.\qed
\endproclaim

The proofs of 3.10 and 3.11 can be found in \cite{FA}, Ch 5.
3.11 is rather elementary, and 3.10 is a fairly easy consequence of
Mori's flip theorem \cite{Mo3}.

\demo{3.12 Proof of 3.3}
Let $X$ be an arbitrary flat closure of $X_K$.
Choose a section $S^\prime \subset X$ of
$X \to S$ (3.7). If $\operatorname{deg} X_\eta \geq 2$, $S^\prime
\subset X$ can be arbitrary, if $\operatorname{deg} X_\eta =1$, we want to
take the base section of the anticanonical system.
Then choose a divisor $B\subset X$ with the following
properties:

3.12.1 $B$ contains $S^\prime$,

3.12.2 $B_\eta \in |-K_{X_\eta}|$ is smooth.

The existence of $B$ follows from elementary properties of Del Pezzo
surfaces.
It is important to understand that, in general, $B \not \in
|-K_X|$.

Choose now a projective desingularization $f: Y \to X$ with the following
properties:

3.12.3 $f_\eta : Y_\eta \overset \cong \to \to X_\eta$ is an isomorphism,

3.12.4 The strict transform $B_Y=f^{-1}_\ast B$ is smooth.

Then the pair $Y, B_Y$ satisfies condition $(\ast)$, and the intent is to run
a minimal model program for $K_Y+B_Y$, and show that the final product
$X^\prime$ is a good model with $-K$ nef. In order to run the minimal model
program, we need to check (3.11) that $B$ is never contracted by a
divisorial contraction. But this is obvious: $K_{Y_\eta}+
{B_Y}_\eta$ is nef, so nothing ever happens to the generic fiber, in
particular $B$ is never contracted.

Let now $X^\prime$ be the end result of the minimal model program. By what
I just said, $K^\prime + B^\prime$ is nef. More than that, since
$(X^\prime_\eta, B^\prime_\eta)=(X_\eta, B_\eta)$, $K^\prime_\eta +
B^\prime_\eta \equiv 0$, which implies that $K^\prime + B^\prime$ is
numerically equivalent to $0$. The subtle point of this
argument, where we use 3.12.1, is to show that actually
$K^\prime +B^\prime \sim 0$ is {\it linearly} equivalent to $0$.
Indeed, by 3.12.1, $B^\prime \to S$ is a relatively minimal smooth
(3.9.2) elliptic surface {\it with a section}. Then $K_{B^\prime}\sim 0$,
which implies $K^\prime +B^\prime \sim 0$. Then $B^\prime \in |-K^\prime |$
is smooth, and $X^\prime$ is a good model with $-K^\prime$ nef.
\qed \enddemo

One can naturally associate, to a scheme $U/S$ with structure
morphism $u:U \rightarrow
S$, and a line bundle $L$ on $U$, a graded
sheaf of ${\Cal O}_S$-algebras ${\Cal R}(U/S,L)= \oplus_{n\geq 0}
u_\ast (nL)$.

If $X$ is a good model with $-K$ nef, it is reasonable to
expect that the sheaf of graded ${\Cal O}$-algebras ${\Cal R}=
\oplus_{n \geq 0} H^0( -nK_X)$ has the usual properties
it enjoys in the theory of Del Pezzo surfaces or Gorenstein elliptic
singularities of surfaces, compare with \cite{R1}.
The next statement shows that this is indeed the case.
The result will be used, in an essential way, in the proof of all our main
results.

\proclaim{3.13 Proposition}
Let $X$ be a good model with $-K$ nef, ${\Cal R}= \oplus_{n \geq 0}H^0
(-nK_X)$ and ${\Cal R}_n = H^0(-nK_X)$ the part of degree $n$.
Let $U$ be the anticanonical model of $X$. Then:

3.13.1 If $d \geq 3$, ${\Cal R}$ is generated by ${\Cal R}_1$.
In particular $-K_X$ is free, and $-K_{U}$ very ample.

3.13.2 If $d= 2$, ${\Cal R}$ is generated by ${\Cal R}_1$
and ${\Cal R}_2$ with one relation. More precisely
${\Cal R}\cong {\Cal O}_C [x, y, z, w]/{(x^2 +f(y,z,w))}$,
where $x$ has degree $2$, $y,z,w$ have degree $1$, and $f$ is homogeneous
of degree $4$. In particular:

3.13.2.1 $-2K_X$ is free and $-2K_U$ very ample,

3.13.2.2  $-K_X$ and $-K_U$ are free,
giving a degree two covering $\pi : U \rightarrow {\Bbb P}^2_S =
\operatorname{Proj} {\Cal O} [y,z,w]$.

3.13.3 If $d=1$, ${\Cal R}$ is generated by ${\Cal R}_1$,
${\Cal R}_2$ and ${\Cal R}_3$ with one relation.
More precisely ${\Cal R} \cong {\Cal O}[x,y,z,w]/{(x^2 +y^3 +yf(z,w)+g(z,w))}$
where $x$ has degree $3$, $y$ has degree $2$, $z$ and $w$ have degree $1$,
$f$ is homogeneous of degree $4$ and $g$ homogeneous of degree $6$.
In particular:

3.13.3.1 $-3K_X$ is free and $-3K_U$ very ample.

3.13.3.2 $-2K_X$ and $-2K_U$ are free,
giving a degree two covering $\pi :U \rightarrow {\Bbb P}^2_S
(2,1,1)$, where ${\Bbb P}^2_S (2,1,1)=\operatorname{Proj}
{\Cal O} [y,z,w]$.

3.13.3.3 The base locus of $|-K_X|$, resp. $-K_U$, consists of a single
section $S_X$, resp. $S_U$ of $X \to S$, resp. $U \to
D$.
Clearly then $X$, resp. $U$ is smooth in a neighborhood of
$S_X$, resp. $S_U$. Let $\epsilon : \tilde U \rightarrow U$ be
the blow up of $S_U$, then $-K_U$
induces a morphism $\tau : \tilde U \rightarrow {\Bbb P}^1_S=
\operatorname{Proj} {\Cal O}_C[z,w]$. A general fiber of $\tau$ is an
elliptic curve.
\endproclaim

3.13 will be a consequence of the following lemma, which is
theorem 3.5 of \cite{R1}. See also \cite{La},
especially the discussion on pages 1270--1278.
It is the key to the classification of Gorenstein elliptic
singularities of surfaces. Recall that a proper effective curve $E$ on a
surface $S$ is said to be {\it numerically m-connected} if for every
decomposition $E=E_1+E_2$ with $E_i>0$, the intersection product
$E_1 \cdot E_2 \geq m$.

\proclaim{3.14 Lemma}
Let $E$ be a numerically 2-connected divisor on
a smooth surface $S$, with $\omega_E = {\Cal O}_E$. Let $L$ be a
line bundle on $E$, such that $\operatorname{deg} L_D \geq 0$ for
each component $D$ of $E$, and let $d=\operatorname{deg} L$.
Set $R= \oplus R_n$ with $R_n= H^0(nL)$.
Then:

3.14.1 If $d \geq 3$, $R$ is generated by $R_1$.
In particular $L$ is generated by global sections.

3.14.2 If $d =2$, $R$ is generated by $R_1$ and $R_2$ with one
relation, and $L$ is generated by global sections.

3.14.3 If $d= 1$, $R$ is generated by $R_1$, $R_2$
and $R_3$ with one relation. $h^0(L)=1$ and the global
sections of $L$
vanish at a single point $P \in E$, lying on a reduced component of $E$,
outside which $L$ is generated by global sections.
\qed \endproclaim

\demo{3.15. Proof of 3.13}
Since $-K_X$ is nef and big, by the vanishing theorem of Kawamata and
Viehweg (\cite{KMM87},
1-2-5 and 1-2-6), $H^1(-nK)=0$ for all $n \geq 0$, so we get exact
sequences:
$$0\rightarrow H^0 \bigl(-nK\bigr)\overset \cdot B\to \to
H^0 \bigl(-(n+1)K \bigr)
\rightarrow H^0 \bigl(-(n+1)K |_B\bigr) \rightarrow 0,$$
combining into a short exact sequence of graded ${\Cal O}$-modules:
$$0 \rightarrow {\Cal R}(X, -K) \overset {\cdot B}\to \rightarrow
{\Cal R}(X, -K) \rightarrow {\Cal R}(B, -K|_B)
\rightarrow 0,$$
where $\cdot B$ is a graded module homomorphism of degree one.
So the statements about $X,-K$ and ${\Cal R}$ follow from the
corresponding statements about $B,L=K|_B$ and ${\Cal R}
(B, L)$.

By definition of good models with $-K$ nef, a general
member $B\in |-K|$ has Du Val singularities,
and $K_B \sim 0$ is linearly equivalent to
zero (adjunction formula). Let $f: B^\prime \to B$ be the minimal
desingularization, then $K_{B^\prime}\sim f^\ast K_B \sim 0$, so $B^\prime_0$
is not a multiple fiber.  Let $E=B^\prime_0$.
Since $B$ is normal, ${\Cal R}({ B}^\prime,L^\prime= f^\ast L)
\cong {\Cal R}(B,L)$, and, using the vanishing theorem as before,
we get exact sequences:
$$0 \rightarrow L^\prime \overset {\cdot E} \to \rightarrow
L^\prime \rightarrow L^\prime|_E \rightarrow 0,$$
combining into a short exact sequence of graded ${\Cal O}$-modules:
$$0 \rightarrow {\Cal R}(B^\prime,L^\prime) \overset {\cdot E}
\to \rightarrow {\Cal R}(B^\prime,L^\prime)
\rightarrow {\Cal R}(E,L^\prime|_E) \rightarrow 0,$$
with the difference that, now, $\cdot E$ is a graded homomorphism of degree
zero.
Everything now follows from 3.14, applied to $E=B^\prime_0$, $L=L^\prime|_E$.
We need to check that the assumptions of 3.14 are satisfied. The conditions
on $L$ are obviously met.
Clearly, $K_E=K_{B^\prime}+E|_E\sim 0+0$, so it
only remains to check that $E$ is numerically 2-connected (this is well
known, but I will explain it anyway).
If $E$ is irreducible, since it is reduced, it is
certainly numerically 2-connected. Otherwise, since the intersection
form is on $E$ is even (every component of $E$ is a $-2$ curve) and
$E$ is not a multiple fiber, $E$ is numerically $2$-connected.
\qed \enddemo

Before we prove theorem 1.10, we need one more result, on the structure of
the cone of curves of a good model with $-K$ nef:

\proclaim{3.16 Cone theorem for good models with $-K$ nef}
Let $X$ be a good model with $-K$ nef, and $\varphi :X \rightarrow U$ the
morphism to the anticanonical model. Then:

3.16.1 $NE(X)\subset N_1(X, {\Bbb R})$ is a
finitely generated rational polyhedral cone (in particular it is closed).

3.16.2 Let  $R \subset NE(X)$ be an extremal ray (in the sense of convex
geometry). Then either
$K_X \cdot R <0$, or $K_X \cdot R =0$ and $E\cdot R <0$ for some
irreducible divisor $E \subset X$.
In particular there is a contraction morphism $\varphi_R :X \rightarrow Y$
associated to $R$. If $\varphi_R$ is small, it is a flopping contraction.

3.16.3 Assume that $X_0$ is not irreducible, and let $E\subset X_0$ be a
component such that $\varphi (E) \subset U$ is a surface. Then, there
is a finite sequence of $E$-flops $X \dasharrow X^\prime$ such that
the proper transform
$E^\prime \subset X^\prime$ is the exceptional
divisor of the contraction of an extremal ray $R \subset NE(X^\prime)$ with
$K_{X^\prime} \cdot R < 0$.
\endproclaim

\demo{Proof}
The result follows formally from the standard cone theorem \cite{KMM87} 4-2-1,
the contraction theorem \cite{KMM87} 3-2-1, and the technical lemma 3.17.
Let me go through this in more detail.

For a divisor $D \subset X$, write $D_{<0} \subset
N_1(X)$ for the halfspace $\{ \,z \mid D \cdot z <0 \,\}$.

Let $E\subset X$ be any effective divisor. Then $K_X+
\epsilon E$ is log terminal for all $0 \leq \epsilon \ll 1$, and
by the cone and contraction theorems (\cite{KMM87} 4-2-1 and 3-2-1)
extremal rays of $ \overline {NE}(X)$
are discrete in the halfspace $(K_X + \epsilon E)_{<0}$ and can be
contracted.

Let $E_i \subset X_0$ be the irreducible components of the central fiber.
In 3.17, it is shown that there are irreducible divisors
$D_j \subset X$ (see 3.17 for the construction of $D_j$) such that, for all
$0< \epsilon \ll 1$:
$$\overline {NE}(X) \subset K_{X<0} \cup
\bigl(\cup_i (K_X+\epsilon E_i)_{<0}\bigr)\cup
\bigl(\cup_j (K_X+\epsilon D_j)_{<0}\bigr)$$
and 3.16.1 follows from the usual compactness arguments.

That extremal rays can be contracted follows from the contraction theorem.
Let $R \subset {NE}(X)$ be an extremal ray, and assume that
$\varphi_R$ is small. Since $-K_X$ is $\varphi_R$-nef, $R$ flips or flops.
But there are no flips on a ${\Bbb Q}$-factorial variety with terminal
singularities of index 1 \cite{Be}. This finishes 3.16.2.

I now prove 3.16.3. Since $X_0$ is reducible,
$E\cdot C <0$ for some curve $C \subset X_0$ and there is an extremal
ray $R$ with $(K_X + \epsilon E)\cdot R <0$. Since $\varphi (E)$ is a surface,
$-K_X |_E$ is nef and big, so either $R$ flops, or $K_X \cdot R <0$,
in which case, by 3.16.2, the
associated contraction is a divisorial contraction contracting $E$. After
a finite sequence of $E$-flops, $E$ can be contracted.
\qed \enddemo

\proclaim{3.17 Lemma}
With the same notation and assumptions as 3.16, the generic fiber
$X_\eta$ is a smooth Del Pezzo surface over $K=K({\Cal O})$. In particular
${NE}(X_\eta)\subset N^1(X_\eta, {\Bbb R})$
is a finitely generated rational polyhedral cone.
Let $R_j$ be the extremal rays,
and $D_{\eta ,j}\subset X_\eta$ irreducible generators of $R_j$.
Also, let $D_j \subset X$ be the closure of $D_{\eta, j}$ in $X$
($D_j$ is an irreducible divisor), and $E_i\subset X_0$ the irreducible
components. If $z$ is a non zero class in $\overline {NE} (X)$, then
either $K_X \cdot z <0$; or $K_X \cdot z=0$ and either $E_i\cdot z<0$
for some $i$, or $E_i \cdot z=0$ for all $i$ and $D_j \cdot z <0$
for some $j$.
\endproclaim

\demo{Proof}
Since $-K_X$ is nef, $K_X \cdot z=0$ if $K_X \cdot z \geq 0$.
Now assume $E_i\cdot z\geq 0$ for all $i$. Let $X_0=\sum \mu_i E_i$
be the central fiber. Since $X_0 \cdot z=0$ and $\mu_i>0$,
this means that $E_i \cdot z=0$ for all $i$.
Note that the anticanonical class $-K_{X_\eta }$ lies
in the interior of the cone $NE(X_\eta)$, so
$-K_{X_\eta}=\sum \alpha_j D_{\eta, j}$
with all $\alpha_j>0$. Because $K_X \cdot z=0$ and $E_i\cdot z=0$ for all $i$,
I have $\sum \alpha_j D_j\cdot z=0$. So, if $D_j\cdot z \geq 0$ for
all $j$, $D_j\cdot z=0$ for all $j$. This is
impossible if $X$ is projective, since clearly the $E_i$s and
$D_j$s generate $N^1(X)$.
\qed \enddemo

\demo{3.18. Proof of 1.10}
Let $X_K$ be a smooth Del Pezzo surface of degree $d\geq 3$, defined over
$K=\operatorname{Frac} {\Cal O}$. We wish to construct a standard model
of $X_K$.
We do so first when $X_K$ is a cubic surface, the general case will follow
from the $d=3$ case.

Assume $d=3$. Choose an arbitrary closure of $X_K$ to a normal variety $X$
over ${\Cal O}$. By 3.4, $X$ is birational to a Gorenstein anticanonical
model $U$. By 3.13.1, the anticanonical system $-K_U$ is very ample, and
defines an embedding $U \hookrightarrow {\Bbb P}^3_{\Cal O}$. Starting with
$U\subset {\Bbb P}^3$, apply the flowchart 2.13. $U$ has {\it canonical}
singularities, so, by 2.20, termination holds, and, by  2.15
($k={\Bbb C}$ so no exceptional model can exist), the program produces
a standard model $X^\prime$. This concludes the proof for $d=3$.

Let now $d\geq 3$ be arbitrary. We will prove the result by descending
induction on $d$. The case $d=3$ being already established, let us assume
$d\geq 4$.
We know (3.7.2) that the
$K$-rational points are dense in the Zariski topology of $X_K$.
Choose a general $K$-rational point $p \in X_K$, and
let $f: Y_K \to X_K$ be the blow up of $p$. It is well
known, and at any rate not difficult to prove, that $Y_K$ is a Del Pezzo
surface of degree $d-1$ over $K$. By the induction assumption, there is a
standard model $Y$ for $Y_K$.
Let $E_{K}$ be the exceptional divisor over $p$, and let $E\subset Y$
be its closure in $Y$. We wish to contract $E$, in
order to get a standard model for $X$.
Let $Y^\prime \to Y$ be a maximal
crepant blow up, $E^\prime$ the proper transform. Since $Y$ has
terminal singularities, $Y^\prime \to Y$ is a ${\Bbb Q}$-factorialization,
and it is a good model with $-K$ nef.
Clearly $E$ is not nef, indeed, on the generic fiber, $E_{K}$ is not
nef. By the cone theorem for good models with $-K$ nef
(3.16), the cone ${NE}(Y^\prime)$ is finitely
generated rational polyhedral, so
there is an extremal ray $R$ with $E \cdot R<0$.
If $K_{Y^\prime} \cdot R <0$, the
contraction $\varphi_R: Y^\prime \to X^\prime$ of $R$ is
divisorial (there are no flipping contractions on a
${\Bbb Q}$-factorial 3-fold with terminal singularities of index 1,
see \cite{Be}) with
exceptional divisor $E$. Otherwise, since $-K_{Y^\prime}$ is nef,
$-K_{Y^\prime}\cdot R=0$ and $R$ is an $E$-flopping ray. There is no
infinite sequence of $E$-flops, therefore, after a finite chain of flops
$Y^\prime \dasharrow Y^{\prime \prime}$, we meet a
divisorial contraction $\varphi: Y^{\prime \prime} \to X^\prime$ as before.
On the generic fiber, this has the effect of contracting the -1 curve
$E_{K}$. We will now check that $X^\prime$ is again a good model with $-K$
nef. This will conclude the proof, since the anticanonical model (3.4) $X$
of $X^\prime$ is then a standard model.
Clearly, $K_{Y^{\prime \prime}}= \varphi^\ast K_{X^\prime}+E$,
which implies that $X^\prime$ has terminal singularities of index 1. Let $C
\subset X^\prime$ be a proper curve.
$\varphi (E)$ is an irreducible curve which is not proper (it maps
isomorphically to $S$), so $C \not \subset \varphi E$ and,
as a cycle, $C = \varphi_\ast C^\prime$, for a curve $C^\prime \subset
Y^{\prime \prime}$, no component of which is contained in $E$.
By the projection formula
then:
$$-K_{X^\prime} \cdot C=\varphi^\ast
(-K_{X^\prime}) \cdot C^\prime=
(-K_{Y^{\prime \prime}}+E)\cdot C^\prime \geq 0,$$
therefore, $X^\prime$ is a good model with $-K$ nef. \qed \enddemo

\newpage

\head 4. Del Pezzo surfaces of degree $2$ \endhead

In this section we construct geometric standard models of Del
Pezzo surfaces of degree 2, i.e., we prove theorem 1.15.
1.15 is an immediate consequence of 3.4, 3.13.2 and 4.1 below
(descending induction on $e, n$).

The material is organized as follows. After stating 4.1 and briefly
discussing its proof, we
begin with some preliminaries on the weighted projective space ${\Bbb
P}(2,1,1,1)$ (4.2), and quartic hypersurfaces in it (4.3). 4.3 is used to
divide up the proof of 4.1 in three main cases 4.5.1--3, which are dealt in
lemmas 4.6--8. The main bulk of this section is devoted to the proof of
4.6--8, using the techniques of \S 2.

The basic set up will be as in \S 3: ${\Cal O}={\Cal O}_{C,p}$
is the local ring at $p\in C$ of a smooth complex curve $C$. As in \S 2,
unless otherwise indicated, all varieties, schemes, morphisms, etc., are
assumed to be defined over $S= \operatorname{Spec} {\Cal O}$. If $X$
is a scheme, $X_\eta$, $X_0$ denote the generic and central fiber.
We only allow birational maps which induce isomorphisms of generic fibers.

We fix a Del Pezzo surface $X_K$, of degree 2 and smooth over $K$, and wish
to construct a standard model $X \subset {\Bbb P}={\Bbb P}(2,1,1,1)_{\Cal
O}$ as in definition 1.13. By 3.4 and 3.13.2, there is a model $U \subset
{\Bbb P}$ with (Gorenstein) canonical singularities.
As in \S 2, let $e(U)$
be the number of crepant exceptional valuations, and $n(U)$ the number of
irreducible components of the central fiber $U_0$.
Our main result in this section is the following:

\proclaim{4.1 Theorem} Let $U \subset {\Bbb P}(2, 1, 1 ,1)_{\Cal O}$ be a
model of $X_K$, with (not necessarily Gorenstein) canonical
singularities. Assume that $U$ is not already a standard model.
Then, one of the following 4.1.1 or 4.1.2 holds:

4.1.1 There exists a birational transformation
$\Phi : {\Bbb P}(2, 1, 1 ,1)_{\Cal O}
\dasharrow {\Bbb P}(2, 1, 1 ,1)_{\Cal O}$
such that, upon setting $U^+=\Phi_\ast U$, $U^+$ has canonical
singularities, $e(U^+) \leq e(U)$ and either:

4.1.1.1 $e(U^+)< e(U)$, or:

4.1.1.2 $e(U^+)=e(U)$ and $n(U^+) < n(U)$.

4.1.2 There exists a birational transformation $U \dasharrow U^+$, where $U^+$
is a standard model.
\endproclaim

The natural approach to constructing standard models in $d=2$ is the same
as for $d=3$: use birational transformations $\Phi :
{\Bbb P}(2, 1, 1 ,1)_{\Cal O}\dasharrow {\Bbb P}(2, 1, 1 ,1)_{\Cal O}$,
centered at the bad singular points of $U$, thus improving the
singularities. The starting point is a model $U$ with {\it Gorenstein}
singularities and, at first, one might hope to preserve this property all
the way through the process. Unfortunately, if the central fiber $U_0$ is
reducible, one is usually forced to introduce some singularities of
index 2 (see
4.6.1), and one is then quickly led to consider arbitrary models $U \subset
{\Bbb P}(2, 1, 1 ,1)_{\Cal O}$. As a further complication, it is sometimes
necessary to even give up ${\Bbb P}(2, 1, 1 ,1)_{\Cal O}$ as an ambient
space, and work in ${\Bbb P}^6_{\Cal O}$ (4.6.3 and example 4.11),
which explains the
awkward formulation of 4.1 into 4.1.1 and 4.1.2. Fortunately, this last
phenomenon only occurs under rather special circumstances, making it
possible to keep it in check. There are many cases to analyze, and,
while trying to give a clear general picture, I will often be quite
sketchy in the arguments.

Before I start the proof of 4.1, I will need to recall some facts about the
weighted projective space ${\Bbb P}(2,1,1,1)$ (4.2) and quartic
hypersurfaces in it (4.3).

\definition{4.2 The space ${\Bbb P}={\Bbb P}(2,1,1,1)$} By definition,
$${\Bbb P}(2,1,1,1)= \operatorname{Proj}k[u,x_1,x_2, x_3],$$
where $u,x_1,x_2, x_3$ have weights, respectively, $2,1,1,1$. Also
${\Bbb P}(2,1,1,1)={\Bbb A}^4_k\setminus \{0\}/k^\times$, where $k^\times$ acts
diagonally on ${\Bbb A}_k^4$,
with coordinate functions $u,x_1,x_2, x_3$, with weights
$2,1,1,1$. As for ordinary projective space, ${\Bbb P}$ is covered by 4
affine charts, corresponding to nonzero values of the coordinate functions.
The important difference is that the point $(1, 0, 0, 0)\in {\Bbb P}$ is
singular. Indeed, the affine chart corresponding to $u\not = 0$ is isomorphic
to the affine toric 3-space ${1\over 2}(1,1,1)$, with inhomogeneous
coordinates which, abusing notation, we still denote $x_1,x_2, x_3$.
The charts with $x_i \not = 0$ are smooth.

It may be useful, for the purpose of drawing pictures and being
able to rely on some geometric
intuition, to realize ${\Bbb P}(2, 1, 1, 1)$ as a cone over the Veronese
surface, via the natural embedding of ${\Bbb P}(2, 1, 1, 1)$ in
${\Bbb P}^6$ with the linear system $|{\Cal O}(2)|$.

The reader can show as an exercise that the automorphism group of ${\Bbb
P}$ is a linear algebraic group of dimension 15, as for ordinary projective
space.

Let ${\Bbb P}_{\Cal O}={\Bbb P}(2,1,1,1)_{\Cal O}$ be the weighted
projective space over ${\Cal O}$. As in \S 2, we will use birational
transformations $\Phi : {\Bbb P}_{\Cal O}
\dasharrow {\Bbb P}_{\Cal O}$, given in coordinates by
$$\Phi^{-1 \ast} (u,x_1, x_2,x_3) = (t^\beta u, t^{\alpha_1} x_1,
t^{\alpha_2} x_2, t^{\alpha_3} x_3).$$
We will say that $(\beta, \alpha_1, \alpha_2, \alpha_3)$ are the {\it
weights} of $\Phi$.
\enddefinition

Let $S$ be a Gorenstein Del Pezzo surface of degree $K_S^2=2$. Due to
classical as well as more recent results, see for instance \cite{R6}, the
anticanonical algebra $\bigoplus_{n\geq 0} H^0(-nK_S)$ defines an embedding
$S \subset {\Bbb P}(2,1,1,1)$ as a hypersurface $(F(u, x_1, x_2, x_3)=0)$
of degree 4. Since $S$ is Gorenstein, necessarily $S$
avoids the singular point $(1, 0, 0, 0) \in {\Bbb P}$, and then, in
suitable coordinates, $F=u^2+G(x_1, x_2, x_3)$, expressing $S$ as the 2-1
cover of ${\Bbb P}^2$ branched along the 4-ic $G=0$. For our purposes, it
will be necessary to examine nonnormal 4-ics $S \subset {\Bbb P}$, which
are not necessarily Gorenstein. The proof of the following result is
entirely elementary:

\proclaim{4.3 Lemma} Let $S=(F=0) \subset {\Bbb P}(2, 1, 1, 1)$
be a nonnormal
hypersurface of degree 4. Then, one of the following is true in suitable
coordinates:

4.3.1 $S$ is reducible or nonreduced and:

4.3.1.1 $x_1$ divides $F$ or $u$ divides $F$ or:

4.3.1.2 $F=Q_1(x_1, x_2, x_3)Q_2(x_1, x_2, x_3)$, where $(Q_i=0)$ are
smooth conics (not necessarily distinct).

4.3.2 $S$ is reduced, irreducible, and:

4.3.2.1 $S$ is Gorenstein and $F=u^2+x_1^2Q(x_1, x_2, x_3)$, and the
singular set of $S$ is the line $u=x_1=0$, or:

4.3.2.2 $K_S$ has index $I=2$, and the
singular set of $S$ contains the line $x_1=x_2=0$.
Here $F=uQ(x_1,x_2,x_3)+G(x_1,x_2,x_3)$, the conic $(Q=0)$ and the quartic
$(G=0)$ have no common component, and the point $(0,0,1)$ is a singular point
on both. \qed \endproclaim

We now begin the proof of 4.1. In coordinates,
$U=(F(u, x_1, x_2, x_3)=0)\subset {\Bbb P}_{\Cal O}$, where
$F \in {\Cal O}[u, x_1, x_2, x_3]$  is (weighted) homogeneous of weight 4.
We adopt the conventions introduced before the proof of 2.23. In particular
we say that $t^k u^\beta x^\alpha \in F$ if the monomial
$u^\beta x^\alpha$ appears in $F$, with coefficient $a \in {\Cal O}$ of
order $m$. $F_0$ is the image of $F$ in $k[u, x_1, x_2, x_3]$.

\definition{4.4 Definition} The {\it axial multiplicity} is the integer
$k=k(U)$, $0\leq k <\infty$, such that $t^ku^2 \in F$.
\enddefinition

It is obvious that the axial multiplicity does not depend of the coordinate
system, and it is therefore a biregular invariant of $U$. Note that $k=0$
if and only if $U$ is Gorenstein.

\definition{4.5 Main division into cases} We prove 4.1 according to the
following division into cases:

{\bf 4.5.1} $U_0$ is reducible or nonreduced.

{\bf 4.5.2} $U_0$ is reduced and irreducible, and $U$ is singular along a
curve.

{\bf 4.5.3} $U_0$ is reduced and irreducible, $U$
has isolated canonical singularities,
but there is a nonterminal point $p \in U$.
\enddefinition

According to 4.5, 4.1 is an immediate consequence of lemmas 4.6--8 below.
In each case the transformation $\Phi : {\Bbb P}_{\Cal O} \dasharrow
{\Bbb P}_{\Cal O}$ is given by specifying its weights. We do our best to
preserve the receptacle ${\Bbb P}_{\Cal O}$ for our spaces. Unfortunately,
this is not always possible, and we wish to call the attention of the
reader on 4.6.3, $k=1$, and example 4.11. Also note 4.6.1, which shows how
nonGorenstein models come into play.

\proclaim{4.6 Lemma} Assume we are in case 4.5.1. In suitable coordinates,
one of the following occurs:

4.6.1 $u|F_0$ and $\Phi : {\Bbb P}_{\Cal O} \dasharrow
{\Bbb P}_{\Cal O}$, with weights $(1, 0, 0, 0)$, satisfies 4.1.1.

4.6.2 $x_1|F_0$ and $\Phi : {\Bbb P}_{\Cal O} \dasharrow
{\Bbb P}_{\Cal O}$, with weights  $(0, 1, 0, 0)$, satisfies 4.1.1.

4.6.3 $F_0=S_1S_2$ as in 4.3.1.2. If $k \geq 2$, $\Phi : {\Bbb P}_{\Cal O}
\dasharrow {\Bbb P}_{\Cal O}$, with weights  $(1 , 1, 1, 1)$, satisfies 4.1.1.
If $k=1$, $U$ has (not necessarily isolated) $cDV$ singularities and there
is $U \dasharrow U^+$ satisfying 4.1.2.
\endproclaim

\proclaim{4.7 Lemma} Assume we are in case 4.5.2. Let $C$ be the singular
set of $U$. Then $U_0$ is singular along $C$, hence it is nonnormal and,
in suitable coordinates, one of the following occurs:

4.7.1 $C=(u=x_1=t=0)$ and $\Phi : {\Bbb P}_{\Cal O} \dasharrow
{\Bbb P}_{\Cal O}$, with weights  $(1, 1, 0, 0)$, satisfies 4.1.1.

4.7.2 $C\supset (x_1=x_2=t=0)$ and $\Phi : {\Bbb P}_{\Cal O} \dasharrow
{\Bbb P}_{\Cal O}$, with weights  $(0, 1, 1, 0)$, satisfies 4.1.1.
\endproclaim

\proclaim{4.8 Lemma} Assume we are in case 4.5.3. In suitable coordinates,
one of the following occurs:

4.8.1 $p=(0,0,0,1) \in {\Bbb P}_0$ and $\Phi : {\Bbb P}_{\Cal O} \dasharrow
{\Bbb P}_{\Cal O}$, with weights  $(2, 1, 1, 0)$, satisfies 4.1.1.

4.8.2 $p=(1,0,0,0) \in {\Bbb P}_0$. In this case either

4.8.2.1 $\Phi : {\Bbb P}_{\Cal O} \dasharrow
{\Bbb P}_{\Cal O}$, with weights  $(1, 1, 1, 1)$, or

4.8.2.2 $\Phi : {\Bbb P}_{\Cal O} \dasharrow
{\Bbb P}_{\Cal O}$, with weights $(0,2,1,1)$,

satisfy 4.1.1.
\endproclaim

In the proof of 4.6--8, we will check that $\Phi$ satisfies 4.1.1,
using 2.21 and the following elementary result:

\proclaim{4.9 Lemma} Let
$U=(F(u, x_1, x_2, x_3,t)=0)\subset {\Bbb P}_{\Cal O}$ be a model of
$X_K$. Let $\Phi : {\Bbb P}_{\Cal O} \dasharrow
{\Bbb P}_{\Cal O}$ be the birational transformation with weights
$(\beta, \alpha_1, \alpha_2, \alpha_3)$, $U^+=\Phi_\ast U$. Assume
that $t^{\beta + \sum \alpha_i}$ divides $F(t^\beta u,
t^{\alpha_1} x_1, t^{\alpha_2} x_2, t^{\alpha_3} x_3)$. Then the
assumptions of 2.21 are satisfied for $\Phi:U \dasharrow U^+$, in other
words:
$$a(\nu_F, K_U)\leq 0$$
for every $\Phi^{-1}$-exceptional divisor $F\subset U^+$.
\endproclaim

\demo{Proof} Let $d$ be maximal such that $t^d$ divides $F(t^\beta u,
t^{\alpha_i}x_i)$. Then $U^+=(F^+=0)$ with
$$F^+(u, x_i,t)={1 \over t^d}F(u, x_i,t).$$
The proof then proceed by comparing differentials on $U$, $U^+$.
For instance, in the chart with $x_3\not =0$ we may use affine coordinates
$v=u/x_3^2, x=x_1/x_3, y=x_2/x_3$, $f(v,x,y,t)=F(v,x,y,1,t)$ and, using
that $F$ is (weighted) homogeneous of degree 4:
$${dvdxdydt \over f(v,x,y,t)} \overset \Phi^{-1 \ast} \to \longrightarrow
{t^{\beta+\sum \alpha_i -4\alpha_3}dvdxdydt \over f(t^{\beta -2\alpha_3} v,
t^{\alpha_1-\alpha_3}x, t^{\alpha_2-\alpha_3}y, t)}=
{t^{\beta+\sum \alpha_i -4\alpha_3}dvdxdydt \over t^{d-4\alpha_3}
f^+(v,x,y,t)}$$
so, if $d \geq \beta + \sum \alpha_i$:
$$K_{U^+}=K_U -\text{effective}$$
as was to be shown. \qed \enddemo

In the proof of 4.6--8, we consider various models $U=(F=0)$. We will
write:
$$F = \sum F_mt^m,$$
where each $F_m=F_m(u, x_1, x_2,x_3)\in {\Cal O}[u, x_1, x_2,x_3]$
is a weighted homogeneous polynomial of degree 4, all of whose coefficients
are units in ${\Cal O}$. In different situations, we will need to put $F$
in suitable normal forms. To simplify the notation, we will usually pretend
that ${\Cal O}$ is a complete local ring. For instance, if the central
fiber $U_0=\bigl(S_1(x_1,x_2,x_3)S_2(x_1,x_2,x_3)=0\bigr)$, with $(S_i=0)$
smooth conics as in 4.3.1.2, and $k=1$, I will say that coordinates can be
chosen so that:
$$F=S_1S_2+u^2t+\sum_{m\geq 1}G_mt^m,$$
where $G_m=G_m(x_1, x_2, x_3)$ are homogeneous of degree 4. Strictly
speaking, this is true only if ${\Cal O}$ is complete. However, the
statement is true in general up to any desired order in $t$, which is
enough for our purposes.

\demo{4.10 Proof of 4.6} By 4.3.1, there are three cases, corresponding to
4.6.1--3:
\medskip
4.10.1 $u|F_0$,

4.10.2 $x_1|F_0$,

4.10.3 $F_0=S_1S_2$ where $(S_i(x_1, x_2, x_3)=0)$ are smooth conics.
\medskip
We discuss each case separately.

In case 4.10.1, in suitable coordinates:
$$F=u(u+Q)+\sum_{m\geq 1} G_m t^m.$$
Then $\Phi$, as given in 4.6.1, obviously satisfies 4.9. Clearly
$k^+=k+1\not =k$, so $U \dasharrow U^+$ is not an isomorphism, and 4.1.1 is
a consequence of 2.21. Note that here $U$ is Gorenstein and $U^+$ is not,
and this is how we are eventually led to consider general models in ${\Bbb
P}_{\Cal O}$, not just those avoiding the vertex.

In case 4.10.2, in suitable coordinates:
$$F=x_1(uL+C)+u^2t^k+\sum_{m\geq 1} (uQ_m+G_m) t^m.$$
Then $\Phi$, as given in 4.6.1, obviously satisfies 4.9. Clearly
$k^+=k-1\not =k$, so $U \dasharrow U^+$ is not an isomorphism, and 4.1.1 is
a consequence of 2.21.

In case 4.10.3, in suitable coordinates:
$$F=S_1S_2+u^2t^k+\sum_{m\geq 1} (uQ_m +G_m) t^m.$$
If $k \geq 2$, $\Phi$, as given in 4.6.3, obviously satisfies 4.9. Clearly
$k^+=k-2\not =k$, so $U \dasharrow U^+$ is not an isomorphism, and 4.1.1 is
a consequence of 2.21.

It remains to discuss the case $k=1$. In suitable coordinates:
$$F=S_1S_2+u^2t+\sum_{m\geq 1}G_mt^m.$$
This is the hardest case, and here we
are forced to give up ${\Bbb P}(2, 1, 1, 1)$ as an ambient space. Let us
show first that $U$ has $cDV$ singularities away from the vertex. First of
all, since $k=1$, the vertex is an isolated quotient singularity of type
${1\over 2}(1,1,1)$ on $U$. Let $p \in U$ be a singular point, other than
the vertex. Changing
coordinates, we may assume that $p=(0,0,0,1)$ and then, in affine
coordinates $v=u/x_3^2, x=x_1/x_3, y=x_2/x_3$, $U=(f=0)$ with:
$$f=s_1(x,y)s_2(x,y)+v^2t+\sum_{m\geq 1}g_m(x,y).$$
Now $s_1, s_2$ are (possibly equal) parameters at $p$, and $v^2t \in f$, so
$p\in U$ is $cA_n$ or $cD_n$ by 2.26. We discuss two cases, the first of
which is the simplest to work with:

4.10.3.1 $S_1\not=S_2$,

4.10.3.2 $S_1=S_2=S$.

In case 4.10.3.1, $U$ has isolated singularities. In other words, the only
thing which is preventing $U$ to be a standard model is that $U_0=T_1+T_2$,
$T_i=(S_i=0)$ is {\it reducible}. Note that each $T_i$ is isomorphic to the
cone $\overline {\Bbb F}_4$ over the rational normal curve of degree $4$.
Let $\varphi:Z \to U$ be a small projective
${\Bbb Q}$-factorialization such that $T_1^\prime$, the proper transform of
$T_1$, is $\varphi$-ample.
Then $T_1^\prime \cong T_1$ is ${\Bbb Q}$-Cartier and it is the
exceptional divisor of an extremal divisorial contraction $\psi: Z \to
U^+$. It is easy to see that $U^+$ is a standard model, i.e., $-2K_{U^+}$
is Cartier and very ample. One way to check this is to construct
$U\dasharrow U^+$ explicitly, as in example 4.11 below, where we also
explain why $U^+$ can not be embedded in ${\Bbb P}(2,1,1,1)$. Another way is to
consider the divisor:
$$D=-2K_Z+T_1^\prime$$
on $Z$. Clearly, $D$ is nef and it restricts to 0 on
$T_1^\prime$. Restricting to a hyperplane section as in the proof of 3.13,
and using 3.14.1, it is easy to see that the ${\Cal O}$-algebra:
$${\Cal R}=\bigoplus H^0(nD)$$
is generated in degree one. In particular $D$ is free from base points and
defines:
$$\psi : Z \to U^+ \subset {\Bbb P}^6_{\Cal O}.$$
Clearly $\psi$ contracts $T_1^\prime$, $U^+$ has terminal singularities,
and $-2K_{U^+}={\Cal O}_{U^+}(1)$ is very ample. In other words, $U^+$ is a
standard model.

In case 4.10.3.2, let $M\geq 0$ be the minimum natural such that $S$ does
not divide $G_M$. Then $U$ has $cD_{M+1}$ singularities generically along
the curve $C=(t=u=S=0)$ (the cases $M=1,2$ are implicitly treated in the
forthcoming discussion).
As in 4.10.3.1, we will explicitly construct a maximal crepant blow up
$\varphi: Z \to U$, and a contraction (of an extremal face) $\psi: Z \to
U^+$, to a standard model $U^+$.

We start off by inductively constructing a partial resolution:
$$Y^{M-1} \to \cdots \to Y^{i+1}\to Y^i \to \cdots \to Y^0=U$$
of the $cD_{M+1}$ curve $C \subset U$, where the exceptional set of each
$Y^{i+1}\to Y^i$ is a smooth irreducible divisor $E^{i+1}$ mapping onto
a curve $C^i \subset Y^i$. Abusing notation we will also denote $E^i$ the
{\it proper} transform of $E^i$ in $Y^j$ for $j\geq i$.

The chain is constructed as follows.

Let $C^0 \subset Y^0=C \subset U$. Let $Y^1 \to Y$ be the blow up of the
ideal sheaf of the central fiber $Y^0_0$ in a small analytic neighborhood
of $C^0$. It is easy to check that the following $(\ast .1.1$--$4)$ hold:

$(\ast .1.1)$ The exceptional set of $Y^1 \to Y^0$ is a smooth
irreducible divisor $E^1
\subset Y^1$, isomorphic to the ruled surface ${\Bbb F}_4$, and $E^1 \to
C^0$ is the ruling on ${\Bbb F}_4$.

$(\ast .1.2)$ Let $E^0 \subset Y^1$ be the proper transform of
the central fiber $Y^0_0
\subset Y^0$. Then $E^1 \cap E^0\subset E^1$ is the negative section on
$E^1 \cong {\Bbb F}_4$.

$(\ast .1.3)$ The singular set of $Y^1$ consists of the vertex and a
section $C^1\subset
E^1$ of the ruling, disjoint from the negative section.

$(\ast .1.4)$ $Y^1$ has $cD_{M}$ singularities generically along $C^1$.

Let $1\leq i <M-2$, and assume that:
$$Y^i\to Y^{i-1}\to \cdots \to Y_0$$
has already been constructed so that the following conditions
$(\ast .i.1$--$4)$ hold:

$(\ast .i.1)$ The exceptional set of $Y^i \to Y^{i-1}$ is a smooth
irreducible divisor $E^i
\subset Y^i$, isomorphic to the ruled surface ${\Bbb F}_4$, and $E^i \to
C^{i-1}$ is the ruling on ${\Bbb F}_4$.

$(\ast .i.2)$ $E^i \cap E^{i-1}\subset E^i$ is the negative section on
$E^i \cong {\Bbb F}_4$.

$(\ast .i.3)$ The singular set of $Y^i$ consists of the vertex and a
section $C^i\subset
E^i$ of the ruling, disjoint from the negative section.

$(\ast .i.4)$ $Y^i$ has $cD_{M+1-i}$ singularities generically along $C_i$.

Then we let $Y^{i+1} \to Y^i$ be the blow up of the ideal sheaf of
$E^i\subset Y^i$. If $i < M-3$, it is easy to check that $(\ast
.i+1.1$--$4)$ hold.
If $i=M-3$, it is easy to check that $(\ast .M-2.1$--$3)$ hold, together with:

$(\ast .M-2.4)$ $Y^{M-2}$ has $cA_3$ singularities generically along
$C_{M-2}$.

Then we let $Y^{M-1} \to Y^{M-2}$ be the blow up of the ideal sheaf of
$E^{M-2}\subset Y^{M-2}$. Once again, $(\ast .M-1.1$--$2)$ hold, together with:

$(\ast .M-1.3)$ The singular set of $Y^{M-1}$ consists of the vertex and a
(possibly reducible) 2-section $C^{M-1}\subset
E^{M-1}$ of the ruling, disjoint from the negative section.

$(\ast .M-1.4)$ $Y^{M-1}$ has $cA_1$ singularities generically along
$C^{M-1}$.

This completes the definition of $Y^{M-1}$. Let, finally, $Y^{M}\to
Y^{M-1}$ be the blow up of the ideal sheaf of $E^{M-1}$, and let $T$ be the
exceptional divisor. If $C^{M-1}$ is irreducible, so is $T$ and we let
$Z=Y^M$. Otherwise $T=T_1+T_2$ and we let $Z \to Y^M$ be a small projective
${\Bbb Q}$-factorialization such that the proper transform of
$T_1$ is relatively ample. Thus, we have constructed our partial
resolution $\varphi : Z \to U$.

It remains to describe the contraction
$\psi: Z \to U^+$.
If $T$ is irreducible, resp. $T=T_1+T_2$ is reducible,
consider the divisor:
$$D=-2K_Z+E^{M-1}+2E^{M-2}+\cdots+ME^0,$$
resp.
$$D=-2K_Z+T_1+2E^{M-1}+\cdots+(M+1)E^0$$
on $Z$. It is easy to check that $D$ is nef and it restricts to 0 on
$E^{M-1}+\cdots +E^0$, resp. $T_1+E^{M-1}+\cdots +E^0$.
Restricting to a hyperplane section as in the proof of 3.13,
and using 3.14.1, it is easy to see that the ${\Cal O}$-algebra:
$${\Cal R}=\bigoplus H^0(nD)$$
is generated in degree one. In particular $D$ is free from base points and
defines:
$$\psi : Z \to U^+ \subset {\Bbb P}^6_{\Cal O}.$$
Clearly $\psi$ contracts $E^{M-1}+\cdots +E^0$, resp. $T_1+E^{M-1}+\cdots
+E^0$, to a point, $U^+$ has terminal singularities,
and $-2K_{U^+}={\Cal O}_{U^+}(1)$ is very ample. In other words, $U^+$ is a
standard model.
\qed \enddemo

\definition{4.11 Example} Let ${\Cal O}={\Bbb C}\{\{t\}\}$, $U=(F=0)\subset
{\Bbb P}_{\Cal O}$ with:
$$F =S_1S_2+t(u^2+G)$$
where $S_1=x_1^2-x_2x_3$, and $S_2=S_2(x_1,x_2,x_3)$, resp.
$G=G(x_1, x_2, x_3)$, are general
homogeneous polynomials of degree 2, resp. 4. This is the situation
described in 4.10.3, $k=1$, $S_1 \not = S_2$.
We will here explicitly construct the
birational transformation $\Phi : U \dasharrow U^+$, to a standard model
$U^+$ (a similar construction is possible in case $S_1=S_2$ as well).

Let $A={\Bbb P}(2,1,1,1)$, and consider the embedding $A \subset {\Bbb
P}^6_{\Cal O}$, with coordinates $u, z_1, ..., z_6$, defined by:
$$\matrix
z_1 &    &z_2 &    &z_3 & &x_1^2-x_2x_3&      &x_1x_2&      &x_2^2 \\
    &z_4 &    &z_5 &    &=&            &x_1x_3&      &x_2x_3&      \\
    &    &z_6 &    &    & &            &      &x_3^2 &      &
\endmatrix$$
$A \subset {\Bbb P}^6_{\Cal O}$ is described by the vanishing of the $2
\times 2$ minors of the matrix:
$$\pmatrix
z_1+z_5 & z_2 & z_4 \\
z_2     & z_3 & z_5 \\
z_4     & z_5 & z_6
\endpmatrix$$
and, in $A$, $U=(F=0)$ with:
$$F=z_1L(z_2, ... z_6)+t(u^2+Q(z_2, ..., z_6))$$
where $L$, resp. $Q$, are homogeneous polynomials of degree 1, resp. 2.

Let $\Phi : {\Bbb P}^6_{\Cal O} \dasharrow {\Bbb P}^6_{\Cal O}$ be the
transformation such that $\Phi^{-1\ast}(u, z_1, z_2, ...,z_6)=
(u, tz_1,z_2, ... z_6)$, $A^+=\Phi_\ast A$, $U^+=\Phi_\ast U \subset A^+$.
$A^+ \subset {\Bbb P}^6_{\Cal O}$ is described by the vanishing of the $2
\times 2$ minors of the matrix:
$$\pmatrix
tz_1+z_5 & z_2 & z_4 \\
z_2     & z_3 & z_5 \\
z_4     & z_5 & z_6
\endpmatrix$$
and $U^+ \subset A^+$ by the equation:
$$F^+=z_1L(z_2, ... z_6)+u^2+Q(z_2, ..., z_6)=0.$$
It is possible, but not easy, to check directly that $U^+$ is a standard
model. The point
here is that $U^+$ lives naturally in $A^+$. Note that the central fiber
$A^+_0$ is a cone (over a cone) over a rational normal curve of degree 4 in
${\Bbb P}^4$, in other words $A^+$ is the total space of a degeneration, in
which the generic fiber $A^+_\eta={\Bbb P}(2,1,1,1)_K$ degenerates to
$A^+_0={\Bbb P}(4,4,1,1)_k$. Because of the structure of its anticanonical
algebra, $U^+_0$ can not be embedded in ${\Bbb P}(2,1,1,1)_k$.
\enddefinition

\demo{4.12 Proof of 4.7}
It is immediate to see that $\Phi$ satisfies 4.9. It is easy to check
that $U \dasharrow U^+$ is not an isomorphism (details are left to the
reader). 4.1.1 is then a consequence of 2.21. \qed \enddemo

4.8 is the analogue of 2.23, and I regret not being able to formulate a
statement that would make this analogy more apparent.
For the proof in case 4.8.2, we will need the
classification of 3-fold terminal singularities \cite{Mo2}, \cite{KSB} \S 6,
see also \cite{R5}. In fact, we only need the index 2 case:

\proclaim{4.13 Theorem} Let $p \in X$ be the (complex analytic) germ of a
3-fold terminal singularity of index 2. Then $p \in X$ is isomorphic to one
of the following germs of hypersurfaces in a suitable affine toric 4-space
$A$:

4.13.1 $xy+ g(z^2, t)=0$ in $A={1 \over 2}(1,1,1,0)$.

4.13.2 $x^2+y^2+ g(z, t)=0$, $g \in m^4$, in $A={1 \over 2}(0,1,1,1)$

4.13.3 $x^2+y^3+ yzt+ g(z, t)=0$, $g \in m^4$, or:

$x^2+y^3+ yzt+ y^n+ g(z, t)=0$, $g \in m^4$ and $n \geq 4$, or:

$x^2+y^3+ yz^2+ y^n+ g(z, t)=0$, $g \in m^4$ and $n \geq 3$, all in
$A={1 \over 2}(1,0,1,1)$.

4.13.4 $x^2+y^3+ yg(z, t)+h(z, t)=0$, $g, h \in m^4$, $h \not = 0$, in $A={1
\over 2}(1,0,1,1)$.
Moreover, each isolated germ 4.13.1--4 is terminal. \qed
\endproclaim

\demo{4.14 Proof of 4.8} In case 4.8.1, choose affine coordinates
$v=u/x_3^2, x=x_1/x_3, y=x_2/x_3, t=t$. Assign weights $2,1,1,1$ to the
variables $v,x,y,t$. Using the technique of the proof of 2.23, it is
possible to show that $p\in U$ has weighted multiplicity $\mu \geq 4$ in
the variables $v,x,y,t$ (details are left to the reader).
It is then immediate that $\Phi$ satisfies 4.9. It is easy to check
that $U \dasharrow U^+$ is not an isomorphism ($k^+\not = k$), and
4.1.1 is then a consequence of 2.21.

I will now discuss 4.8.2 in detail. $U=(F=0)$ with:
$$F= \sum (uQ_m+G_m)t^m+t^ku^2.$$
The singularity at the origin is:
$$\bigl(\sum (Q_m+G_m)t^m+t^k=0\bigr) \subset {1 \over 2}(1,1,1,0)$$
($t$ has weight 0). First of all $k \geq 2$, otherwise $p \in U$ is a
terminal quotient singularity of type ${1 \over 2} (1,1,1)$.
If $Q_0=0$, $\Phi$ as given in 4.8.2.1 satisfies 4.9, $k^+\not = k$, and
4.1.1 is then a consequence of 2.21. We may then assume that $Q_0\not = 0$.
In suitable coordinates then $Q_0=x_1^2$, otherwise $p \in U$ is a terminal
singularity of type 4.13.1. Changing coordinates again, we may assume that
$Q_1=Q_1(x_2,x_3)$ does not depend on $x_1$.
Now $k \geq 3$, otherwise $p \in U$ is a terminal
singularity of type 4.13.2.
If $Q_1\not = 0$, $p \in U$ is
a terminal singularity of type 4.13.3, so we may assume $Q_1=0$. We will
now show that $k \geq 4$. Indeed, if $k=3$, $G_0(0,x_2,x_3)=0$, otherwise
$p \in U$ is a terminal
singularity of type 4.13.4, but then $U_0$ is reducible,
contradicting the assumptions. Then $\Phi$, as given in 4.8.2.2,
satisfies 4.9, $k^+\not = k$, and
4.1.1 is then a consequence of 2.21.\qed \enddemo

\newpage

\head 5. Examples of local rigidity \endhead

In this section, we prove theorem 1.18. 1.18 is an immediate consequence of
5.1, 5.2. 5.1 states that if a Mori fiber space $X \to
\operatorname{Spec} {\Cal O}$ is birational to a Gorenstein standard model
$X^\prime \to \operatorname{Spec} {\Cal O}$, {\it there exists} a surface
$B \equiv -K_X$, such that $B_\eta$ is smooth but $B$ has worse than Du Val
singularities. On the other hand, if $X\to \operatorname{Spec} {\Cal O}$ is
a standard model, by definition the central fiber is reduced, so every
$B\equiv -K_X$ is also linearly equivalent to $-K_X$, and 5.2 states that
if $X \to \operatorname{Spec} {\Cal
O}$ is ``sufficiently general'' of index $I\geq 2$, then
{\it every} member $B \in |-K_X|$, such that $B_\eta$ is smooth, has Du Val
singularities.

The proof of 5.1 uses the N\"other-Fano inequalities, as stated and proved
in \cite{Co}, 4.2. The proof of 5.2 is a tedious explicit calculation in
coordinates which occupies the main bulk of this section (lemmas 5.4--7).

In 5.8, we give an example showing that standard models are not
unique. We work with ${\Cal O}={\Bbb C}\{\{t\}\}$.

\proclaim{5.1 Theorem} Let $X_K$ be a Del Pezzo surface, $X \to
\operatorname{Spec} {\Cal O}$, $X^\prime \to
\operatorname{Spec} {\Cal O}$ be models of $X_K$ and $\Phi : X \dasharrow
X^\prime$:
$$
\definemorphism{birto}\dashed \tip \notip
\diagram
   X \rbirto^\Phi \dto & X^\prime \dto \\
\operatorname{Spec} {\Cal O} & \operatorname{Spec} {\Cal O}\\
\enddiagram
$$
a birational map, which (according to our usual conventions)
restricts to the identity of generic fibers. Assume:

5.1.1 $X \to \operatorname{Spec} {\Cal O}$ is a Mori fibration, and:

5.1.2 $X^\prime \to \operatorname{Spec} {\Cal O}$ is a standard model of
index $I=1$ (i.e., it is a Gorenstein standard model).

Then there is a surface $B \equiv -K_X$, such that $B_\eta$ is smooth but
$K+B$ is not canonical, or in other words, $B$ has worse than Du Val
singularities.
\endproclaim

\demo{Proof} $|-K_{X^\prime}|$ is very ample. The result is then an
immediate consequence of the N\"other-Fano inequalities \cite{Co}, 4.2.
\qed \enddemo

\proclaim{5.2 Theorem} Let $X_{\Cal O}$ be a general standard model of degree
$d=1$
or $d=2$, and of index $I \geq 2$. Let $B \in |-K_X|$ be a member such that
$B_\eta$ is smooth. Then $B$ has Du Val singularities.
\endproclaim

\definition{5.3 Remark} In 5.4--7, we explicitly
spell out the conditions $X$ needs to satisfy, in order for the statement
5.2 to be true.
\enddefinition

Let $X \subset {\Bbb P}_{\Cal O}$ be as in 5.2. If $d=2$, ${\Bbb P}=
{\Bbb P}(2,1,1,1)$ and we choose homogeneous coordinates $(u, x_1, x_2,
x_3)$ with weights $(2,1,1,1)$. Then $X$ is defined by an equation:
$F(u,x_1, x_2, x_3, t)=0$, where as usual we write:
$$F=\sum F_m t^m,$$
with $F_m=F_m(u, x_1, x_2, x_3)$ is (weighted) homogeneous of degree 4.
If $d=1$, ${\Bbb P}={\Bbb P}(3,2,1,1)$ and we choose homogeneous
coordinates $(u,v, x_1, x_2)$ with weights $(3,2,1,1)$, then $X$ is a
sextic and similar remarks apply. 5.2 is a direct consequences of the
following more precise results 5.4--7.

\proclaim{5.4 Lemma} Let $X$ be as in 5.2, with $d=2$, $I=2$. Then in
suitable coordinates
$$F=uQ(x_1, x_2, x_3, t)+G(x_1, x_2, x_3, t)+t^ku^2,$$
with $1 \leq k <\infty$. Assume:

5.4.1 the conic $Q_0(x_1, x_2, x_3) =0$, in ${\Bbb P}^2$ with coordinates
$x_1, x_2, x_3$, is a smooth conic, and intersects transversally the 4-ic
$G_0=0$.

5.4.2 $k=1$.

Let $B \in |-K_X|$ be a member such that
$B_\eta$ is smooth. Then $B$ has Du Val singularities.
\endproclaim

\proclaim{5.5 Lemma} Let $X$ be as in 5.2, with $d=1$, $I=2$. Then in
suitable coordinates
$$F=u^2 +v^2 Q(x_1, x_2, t)+vG(x_1, x_2, t)+H(x_1, x_2, t)+t^kv^3,$$
where $1 \leq k <\infty$ and $Q, G, H$ have degrees $2, 4, 6$. Assume:

5.5.1 $Q_0=0$, $G_0=0$ and $H_0=0$ do not have a common solution
in ${\Bbb P}^1$ with coordinates $x_1, x_2$.

5.5.2 $k=1$.

Let $B \in |-K_X|$ be a member such that
$B_\eta$ is smooth. Then $B$ has Du Val singularities.
\endproclaim

\proclaim{5.6 Lemma} Let $X$ be as in 5.2, with $d=1$, $I=3$.
Then in
suitable coordinates
$$F=uvL(x_1,x_2, t) +uC(x_1, x_2, t)+v^3+ vG(x_1, x_2, t)+H(x_1, x_2,
t)+t^ku^2,$$
where $1 \leq k <\infty$, $1 \leq s <\infty$ and
$L, C, G, H$ have degrees $1,3, 4, 6$. Assume:

5.6.1 $(L_0=C_0=G_0=H_0=0)=\emptyset$ in ${\Bbb P}^1$
with coordinates $x_1, x_2$.

5.6.2 $k=1$.

Let $B \in |-K_X|$ be a member such that
$B_\eta$ is smooth. Then $B$ has Du Val singularities.
\endproclaim

\proclaim{5.7 Lemma} Let $X$ be as in 5.2, with $d=1$, $I=6$. Then in
suitable coordinates
$$\align
F=uvL(x_1,x_2, t) +v^2 Q(x_1, x_2, t) +uC(x_1, x_2, t)+
vG(x_1, x_2, t)&+H(x_1, x_2,t)+\\
               &+t^ku^2+t^su^2,
\endalign$$
where $1 \leq k <\infty$ and $L, Q,C, G, H$ have degrees $1,2,3, 4, 6$. Assume:

5.7.1 $(L_0=Q_0=C_0=G_0=H_0=0)=\emptyset$ in ${\Bbb P}^1$
with coordinates $x_1, x_2$.

5.7.2 $k=s=1$.

Let $B \in |-K_X|$ be a member such that
$B_\eta$ is smooth. Then $B$ has Du Val singularities.
\endproclaim

The proof of 5.4--7 is an elementary explicit calculation.
As an example, we prove 5.4 an 5.5.

\demo{Proof of 5.4} After a change of coordinates (this always preserves
5.4.1, 5.4.2) we may assume that $B=(x_3=0)$ and $Q_0(x_1,x_2,0)=x_1 x_2$
or $x_1^2$.

If $Q_0(x_1,x_2,0)=x_1 x_2$, $B_0=(ux_1x_2+G_0(x_1,x_2,0)=0)$ in ${\Bbb
P}(2,1,1)$ with coordinates $u,x_1,x_2$. By 2.23.1, $B$ has $A_n$
singularities away from $(1,0,0)$ and, by 5.4.2,
it has an $A_1$ singularity there.

If $Q_0(x_1,x_2,0)=x_1^2$, $B_0=(ux_1^2+G_0(x_1,x_2,0)=0)$ in ${\Bbb
P}(2,1,1)$ with coordinates $u,x_1,x_2$. By 5.4.1, $x_1^2$ does not divide
$G_0(x_1,x_2,0)$. Then $B$ is smooth everywhere except possibly at $(1,0,0)$
and, by 5.4.2, it has an $A_1$ singularity there. \qed \enddemo

\demo{Proof of 5.5} After a change of coordinates (this always preserves
5.5.1, 5.5.2) we may assume that $B=(x_2=0)$ and $B_0$ is given by:
$$u^2+v^2Q_0(x_1,0)+vG_0(x_1,0)+H_0(x_1,0)=0$$
in ${\Bbb P}(3,2,1)$ with coordinates $u,v,x_1$. By 5.5.1 and 2.23.1, $B$
has $A_n$ singularities except possibly at $(1,0,0)$ and, by 5.5.2,
it has an $A_1$ singularity there. \qed \enddemo

\definition{5.8 Example} We show that standard models are not unique.
Let $X/{\Cal O}$ be a standard model of degree $d=3$. Assume that the
morphism $X \to \operatorname{Spec} {\Cal O}$ is smooth, and the central
fiber $X_0$ contains an {\it Eckhard point}. In other words, there is a
point $p\in X_0$ and a hyperplane section $p \in B_0=\ell_1+\ell_2+\ell_3
\subset X_0$, consisting of 3 lines $\ell_i$, all passing through $p$.
Note that having an Eckhard point is a codimension one property in the
space of cubic surfaces. I will construct a birational map $\Phi : X
\dasharrow X^\prime$, to a new standard model $X^\prime$.
Let $f: Z \to X$ be the blow up of $p$, $F$ the $F$-exceptional divisor,
$G=f^{-1}_\ast X_0$. The normal bundle, in $Z$, of the transform
$\ell^\prime_i$ of $\ell_i$, is ${\Cal O}(-1)\oplus {\Cal O}(-2)$. The
curves $\ell_i^\prime$ can therefore be simultaneously inverse flipped, let
$Z \dasharrow Z^\prime$ be the inverse flip, $F^\prime$, $G^\prime$ the
proper transforms. It is easy to see that $-K_{Z^\prime} + G^\prime$ is nef
on $Z^\prime$ and defines a morphism $g:Z^\prime \to X^\prime$
contracting $G^\prime$, and $X^\prime $ is a standard model.
$X^\prime_0$ is a normal cubic surface with a $D_4$ singularity.
It must be possible to realize explicitly
$X, X^\prime \subset {\Bbb P}^3_{\Cal O}$,
and $X \dasharrow X^\prime$ from an explicit birational map
${\Bbb P}^3_{\Cal O} \dasharrow {\Bbb P}^3_{\Cal O}$, but my attempts at
this have been unsuccessful.
\enddefinition

\enddocument